\documentclass[11pt]{article}
\usepackage{amsmath}
\usepackage{bbm}
\usepackage{amsfonts}
\usepackage{graphicx}
\usepackage{enumerate,booktabs}
\usepackage{mathtools}
\usepackage[sort,compress]{cite}
\usepackage{parskip}
\usepackage[papersize={8.5in,11in},margin=1in]{geometry}
\usepackage{caption}
\usepackage{subcaption}
\usepackage{xcolor}

\usepackage{hyperref}
\hypersetup{
  pdfcreator = {},
  pdfproducer = {}
}

\newtheorem{theorem}{Theorem}

\newtheorem{example}{Example}

\numberwithin{equation}{section}

\newcommand*{\affaddr}[1]{#1} 
\newcommand*{\affmark}[1][*]{\textsuperscript{#1}}
\newcommand*{\email}[1]{\texttt{#1}}

\begin{document}

	\title{Specification of additional information for solving stochastic inverse problems}

\author{%
Wayne Isaac T. Uy\affmark[1], Mircea D. Grigoriu\affmark[1,]\affmark[2]\\
\affaddr{\affmark[1]Center for Applied Mathematics}\\
\affaddr{\affmark[2]Department of Civil and Environmental Engineering}\\
\email{\{wtu4,mdg12\}@cornell.edu}\\
\affaddr{Cornell University}%
}
\date{\vspace{-3ex}}

\maketitle


\begin{abstract}
  Methods have been developed to identify the probability distribution of a random vector $Z$ from information consisting of its bounded range and the probability density function or moments of a quantity of interest, $Q(Z)$. The mapping from $Z$ to $Q(Z)$ may arise from a stochastic differential equation whose coefficients depend on $Z$. This problem differs from Bayesian inverse problems as the latter is primarily driven by observation noise. We motivate this work by demonstrating that additional information on $Z$ is required to recover its true law. Our objective is to identify what additional information on $Z$ is needed and propose methods to recover the law of $Z$ under such information. These methods employ tools such as Bayes' theorem, principle of maximum entropy, and forward uncertainty quantification to obtain solutions to the inverse problem that are consistent with information on $Z$ and $Q(Z)$.  The additional information on $Z$ may include its moments or its family of distributions. We justify our objective by considering the capabilities of solutions to this inverse problem to predict the probability law of unobserved quantities of interest.
\end{abstract}

\section{Introduction} \label{SISC:sec:Intro}

Inverse problems emerge from applications in science and engineering when information about inputs to a system is sought given measurements of observable quantities. Suppose that the physical system is modeled by the mapping $A(x,Z) \mapsto U(x,Z)$ where $x \in \mathbb{R}^d$ is the spatial variable, $Z \in \mathbb{R}^n$ is a parameter, $A(x,Z)$ is a deterministic function, and $U(x,Z)$ is the response. One of the most commonly investigated aspects of this mapping involves deterministic inverse problems. They deal with the estimation of the unknown parameter $Z$ provided measurements $\{U_{x_i}\}_{i=1}^{N_{obs}}$ of $U$ at spatial points $x_i \in \mathbb{R}^d$ which are referred to as quantities of interest. The ill-posedness of this problem is due to the non-uniqueness of the solution for $Z$ and is commonly addressed via two well-established methods. Optimization approaches \cite{book:Hansen2010} solve for $Z$ by minimizing the objective function
$\sum_{i=1}^{N_{obs}} |U_{x_i} - U(x_i,Z)|^2 + \lambda^2 \|Z\|^2$
where the regularization term $\lambda^2 \|Z\|^2$ suppresses noisy solutions. In contrast, Bayesian approaches \cite{book:KaipioS2005} construct the solution as a probability density function (pdf) for $Z$ instead of a point estimate, i.e. we acquire a probabilistic solution to a deterministic problem, by specifying a prior pdf on $Z$ based on available information on the unknown parameter. The observations $\{U_{x_i}\}_{i=1}^{N_{obs}}$ are typically assumed to be polluted by random noise whose law coupled with the prior on $Z$ yield the posterior pdf on $Z$. Despite the connections between both approaches, they emphasize that additional information on $Z$ is required to address the ill-posedness that arises from solving the inverse problem.

In this work, we focus instead on a different class of inverse problems where the unknown quantity is inherently stochastic. Using the same mapping as above, let $Z$ be a random vector defined on the probability space $(\Omega,\mathcal{F},P)$ such that $A(x,Z)$ is a random field while $U(x,Z)$ is the stochastic response. The stochastic inverse problem we address is the following \cite{paper:BreidtBE2011,paper:ButlerETDW2014,
paper:ButlerGEDW2015,paper:GerbeauLT2018}:
\begin{quote}
	Determine the probability law of $Z$ given probabilistic information (such as pdf) of the quantity of interest $Q(Z) \in \mathbb{R}^m$ represented by functionals of the stochastic response $U(x,Z)$.
\end{quote}
By definition, the quantity of interest $Q$ is a function of the response $U$. Examples include $Q(Z) = \text{max}_x |U(x,Z)|$ and $Q(Z) = (U(x_1,Z),\dots,U(x_m,Z))$ where $x_1,\dots,x_m \in \mathbb{R}^d$ are fixed. Also, the output $Q(Z)$ may or may not contain observation noise; regardless, it is still a random quantity as it depends on $Z$ which distinguishes it from what is encountered in the Bayesian formulation above in which $Q(Z)$, without the additive observation noise, is deterministic. This inverse problem is the direct reverse of forward uncertainty propagation which is concerned with obtaining probabilistic information of $Q(Z)$ given the probability law of $Z$.

Without additional information on $Z$, solving this stochastic inverse problem becomes challenging since distinct probability laws on $Z$ can produce the same law for $Q(Z)$ cf. \cite[pp. 1839-1840]{paper:BreidtBE2011}. This ill-posedness is then compensated by imposing assumptions on the law of $Z$ \cite{paper:BreidtBE2011, paper:ButlerETDW2014, paper:ButlerGEDW2015,paper:BorggaardV2015,paper:ZabarasG2008,paper:NarayananZ2004,paper:GerbeauLT2018}.
We mainly scrutinize two methods proposed in literature in which the only information known about $Z$ comprises, at most, its bounded range. The first approach \cite{paper:BreidtBE2011, paper:ButlerETDW2014, paper:ButlerGEDW2015} uses the disintegration theorem for probability measures to obtain a pdf for $Z$ given the pdf of $Q(Z)$. The second approach \cite{paper:BorggaardV2015} aims to approximate the unknown random field $A(x,\omega)$, $\omega \in \Omega$, by solving an optimization problem; it is shown that this is conceptually identical to the stochastic inverse problem we consider above. 


%
%
%

These existing methodologies are examined and the inadequacy of their basic implementations in recovering the true law of $Z$ in the absence of further information on $Z$ is used to motivate the following contribution of this work:

\begin{quote}
	\emph{Identification of additional information to recover the true law of $Z$.}
\end{quote}
This work offers supplemental analysis to \cite{paper:BreidtBE2011, paper:ButlerETDW2014, paper:ButlerGEDW2015, paper:ButlerJW2018} in that while the same inverse problem is tackled, the desired properties that we seek in designing the solution are different. Our focus here is on the reconstruction of the true law which is made possible by incorporating various types of additional information on $Z$ not considered in \cite{paper:BreidtBE2011, paper:ButlerETDW2014, paper:ButlerGEDW2015, paper:ButlerJW2018}. We investigate what type of information $Z$ needs to be equipped with and the corresponding tools that can be employed to solve the inverse problem consistently provided such additional information. Furthermore, we underscore the importance of recovering the true law of $Z$ so that the resulting distribution on $Z$ can be used to characterize unobserved quantities of interest other than those to which it was calibrated. As such, the contributions of this work aid in improving predictions based on solutions to stochastic inverse problems. The examples analyzed here are not meant to be critiques of existing literature. Instead, they offer further insight to better understand how these methods operate and suggest clues on issues that need to be accommodated to enhance their applicability.


\section{Absence of information on the unknown random quantity}\label{SISC:sec:AbsentInfo}

We survey existing methods based on: the disintegration theorem (Section~\ref{SISC:subsec:BETMethod}) and an optimization approach (Section~\ref{SISC:subsec:VanWyk}), to tackle the inverse problem. The information available on the unknown random quantity is limited to its bounded domain. It is argued that this is insufficient to recover the true law of the unknown quantity. 

\subsection{Disintegration of probability measures on generalized contours}\label{SISC:subsec:BETMethod}

The disintegration theorem is summarized in Section~\ref{SISC:subsubsec:BETreview} which serves as the basis of the works \cite{paper:BreidtBE2011, paper:ButlerETDW2014, paper:ButlerGEDW2015}. Two design approaches for the disintegration theorem are then elaborated in Sections~\ref{SISC:subsubsec:DesignLebesgueMeas} and~\ref{SISC:subsubsec:DesignNonLebesgue} with each containing an example applying the method. The examples underscore that the ansatz imposed by this methodology -- i.e., the pdf on the generalized contour is uniform, may be restrictive as the pdf on the contours can possess complex behavior. They also show that the true pdf of the unknown quantity is under/overestimated if this ansatz is accepted.


\subsubsection{Review of methodology}\label{SISC:subsubsec:BETreview}

For the random vector $Z \in \mathbb{R}^n$ defined on the probability space $(\Omega,\mathcal{F},P)$, let $Q$ be the mapping $Q: \Gamma \rightarrow \mathcal{D}$ where  $\Gamma \coloneqq Z(\Omega)$ is a compact subset, $\mathcal{D} \coloneqq Q(\Gamma) \subset \mathbb{R}^m$ with $m < n$, and $Q(z) = q(U(x,z))$ for $z \in \Gamma$ and for some function $q(\cdot)$ that is locally differentiable. In the rest of this work, distinction is made between pdfs constructed with respect to Lebesgue and non-Lebesgue measures. For measurable $A \subset \Gamma$ and  $B \subset \mathcal{D}$, denote by $P_Z(A) = \displaystyle \int_A \rho_Z(z) \, d \mu_Z$ and $P_Q(B) = \displaystyle \int_B \rho_Q(q) \, d \mu_Q$ the probability measures on $Z$ and $Q$, respectively, where $\rho_Z(z)$ and $\rho_Q(q)$ are the corresponding pdfs with respect to the specified measures  $\mu_Z$ and $\mu_Q$. Unlike in calculus-based probability, $\mu_Z$ and $\mu_Q$ are not restricted to be Lebesgue. For notation purposes, if the pdfs are constructed with respect to the Lebesgue measure, they are denoted by $f$ instead of $\rho$. 



The method proposed in \cite{paper:BreidtBE2011, paper:ButlerETDW2014, paper:ButlerGEDW2015} then addresses the following:
\begin{quote}
	Given the pdf of $Q(Z)$ and the bounded range $\Gamma$ of $Z$, estimate the pdf of $Z$.
\end{quote}
The cited approaches compute probabilities in $\Gamma$ by viewing the domain from a coordinate system of contours instead of the traditional Cartesian system. A summary of \cite{paper:ButlerETDW2014} is outlined below.

\begin{figure}[!h]
\centering
\begin{subfigure}{.5\textwidth}
  \centering
  \includegraphics[width=.8\linewidth]{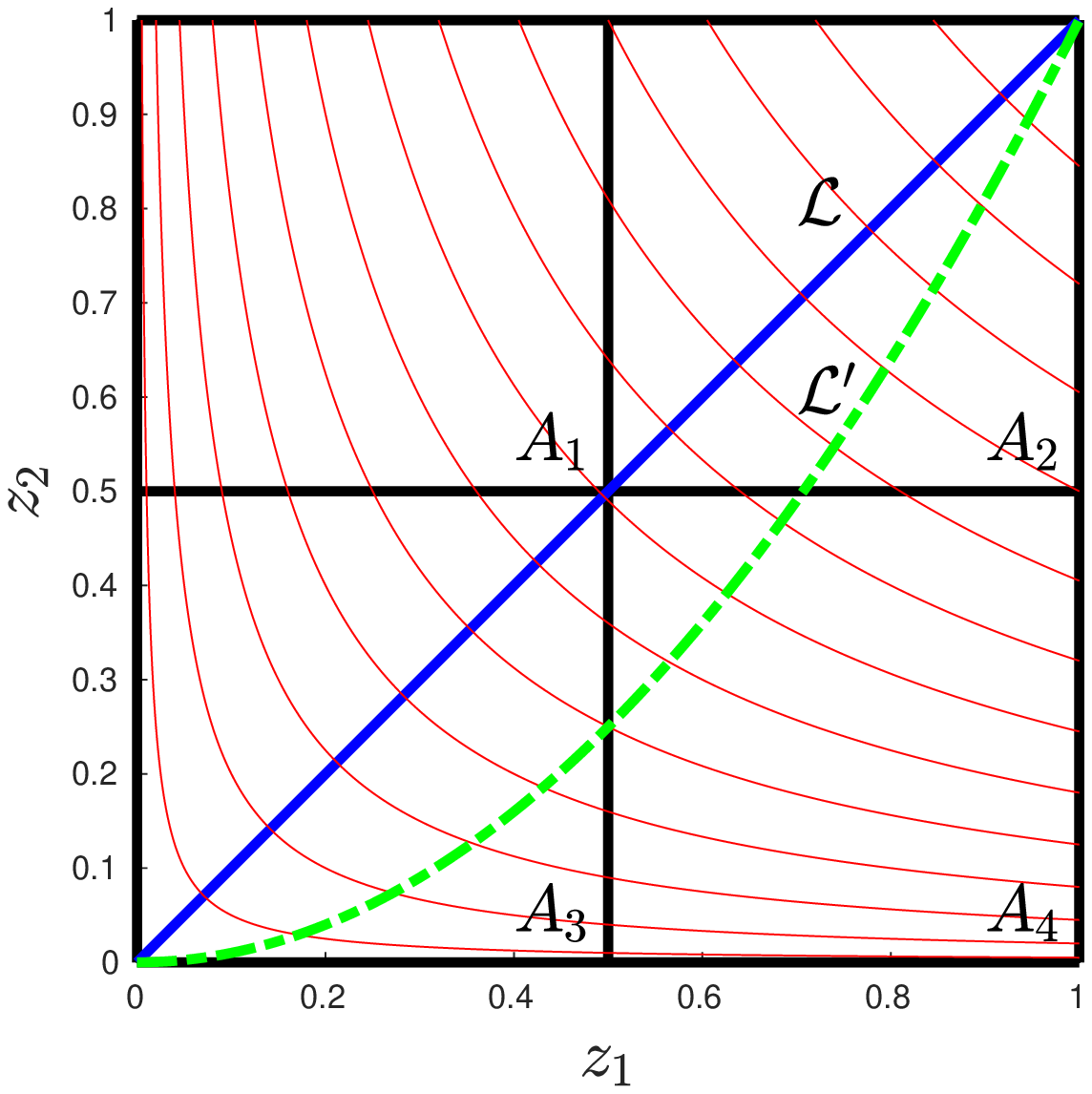}
  \caption{}
  \label{SISC:fig:new_coord_system}
\end{subfigure}%
\begin{subfigure}{.5\textwidth}
  \centering
  \includegraphics[width=0.8\linewidth]{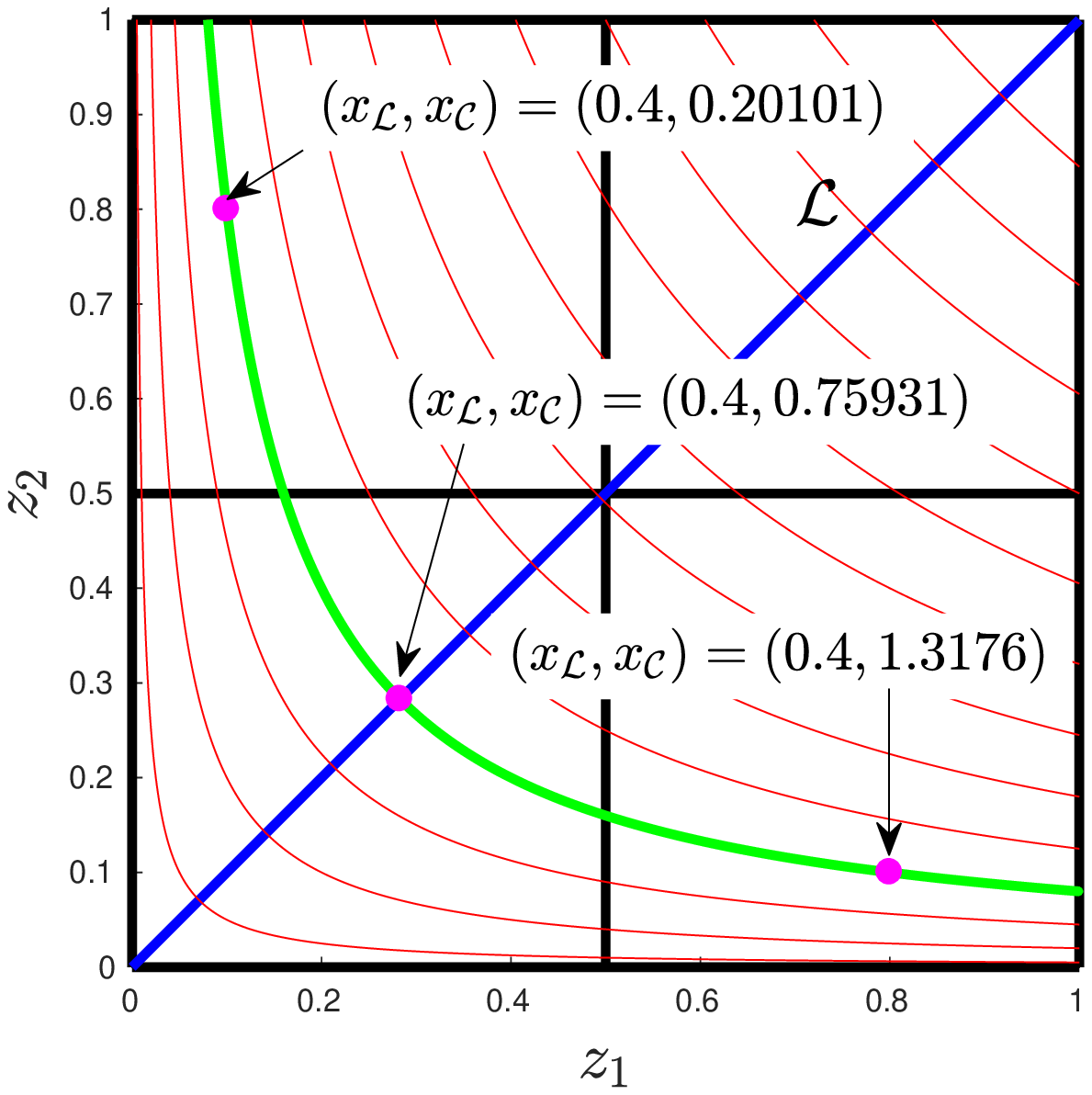}
   \caption{}
   \label{SISC:fig:new_coord_system_w_examples}
\end{subfigure}
\caption{Left panel: Generalized contours (thin red lines) and the two transverse parameterizations (thick blue and dashed green line) for $Q(z_1,z_2)=z_1 \cdot z_2$. Right panel: Illustration of the change in coordinate system from $z\in \Gamma$ to $(x_{\mathcal{L}},x_{\mathcal{C}})$.}
\label{SISC:fig:new_coord_main}
\end{figure}

For $d \in \mathcal{D}$, define the set of points $\{z \in \Gamma | Q(z) = d\}$ to be a generalized contour. This is a generalization of contour curves when $n=2, m = 1$. Generalized contours are equivalence classes with the relation $a \sim b$ if and only if $Q(a) = Q(b)$. Consequently, as $\mathcal{D}$ is the range of $Q$, the domain $\Gamma$ is a union of generalized contours. A representative element from each generalized contour can then be selected which serves as an indexing mechanism across all contours. This indexing set is a $m$-dimensional manifold that intersects each contour once and is called a transverse parameterization. In other words, there is a bijection between points in the transverse parameterization and the generalized contours. To clarify these concepts, an example is shown in Figure~\ref{SISC:fig:new_coord_system} below for  $n = 2, m = 1$ and $Q(z_1,z_2) = z_1 \cdot z_2$. The red thin curves are the contours of $Q$ while the solid blue curve marked $\mathcal{L}$ ($z_2 = z_1$) and the dashed green curve marked $\mathcal{L}'$ ($z_2 = z_1^2$) are different transverse parameterizations. Because the transverse parameterization is not unique in general, we will only consider one parameterization in what follows and denote it by $\mathcal{L}$. 

Every $z \in \Gamma$ can then be expressed under this new coordinate system. Let $\pi: \Gamma \rightarrow \mathcal{L}$ be the onto mapping such that for $x_{\mathcal{L}} \in \mathcal{L}$, $\pi^{-1}(x_{\mathcal{L}})$ is the corresponding generalized contour and that for $A \subset \Gamma$, $\pi(A)$ is the portion of the transverse parameterization $\mathcal{L}$ which intersects the generalized contours contained in $A$. We associate every $z \in \Gamma$ with $(x_{\mathcal{L}},x_{\mathcal{C}})$ where $x_{\mathcal{L}} \in \mathcal{L}$ represents the generalized contour in which $z$ resides and $x_{\mathcal{C}}$ is the coordinate along the generalized contour $\pi^{-1}(x_{\mathcal{L}})$. Figure~\ref{SISC:fig:new_coord_system_w_examples} exhibits this change of coordinate system where for this example, $x_{\mathcal{L}}$ parametrizes the arc length of the transverse from the origin while $x_{\mathcal{C}}$ parametrizes the arc length of the contour from $z_2 = 1$. The transverse $\mathcal{L}$ is parameterized by $z_2 = z_1$. The green thick curve is the contour $\pi^{-1}(x_{\mathcal{L}})$ that is indexed by $x_{\mathcal{L}} = 0.4$, that is, $x_{\mathcal{L}} \in \mathcal{L}$ is 0.4 units from the origin. Meanwhile, the $x_{\mathcal{C}}$-coordinate of each of the 3 magenta circles on this contour is obtained by measuring the arc length of said contour from $z_2 = 1$ up to each of the 3 marked circles. 

With this new coordinate system, the randomness in $Z$ implies that there are random vectors $X_{\mathcal{L}}$ and $X_{\mathcal{C}}$ associated with $x_{\mathcal{L}}$ and $x_{\mathcal{C}}$, respectively. In particular, for measurable $K \subset \mathcal{L}$, denote by $P_{X_{\mathcal{L}}}(K) = \displaystyle \int_K \rho_{X_{\mathcal{L}}} (x_{\mathcal{L}}) \, d\mu_{X_{\mathcal{L}}}$ the probability measure on $X_{\mathcal{L}}$ where $\rho_{X_{\mathcal{L}}}$ is the pdf on $\mathcal{L}$ with respect to $\mu_{X_{\mathcal{L}}}$ that is specified. To solve the inverse problem, the probability of a measurable set $A \subset \Gamma$ can therefore be computed using the disintegration theorem for probability measures as follows \cite[Corollary 4.1, Theorems 4.4, 4.5]{paper:ButlerETDW2014}.

\begin{theorem}[Disintegration theorem] \label{SISC:thm:DisintThm}
Let $P_Z$ be the probability measure defined by the law on $Z$ and $Q$ the measurable mapping between $\Gamma$ and $\mathcal{D}$. For measurable $A \subset \Gamma$, $P_Z$ admits the following disintegration
\begin{align}\label{SISC:eq:disintegratn_thm}
P_Z(A) = \displaystyle \int_{\pi(A)} \int_{\pi^{-1}(x_{\mathcal{L}})\cap A} \rho_{X_{\mathcal{C}} | X_{\mathcal{L}}} (x_{\mathcal{C}} | x_{\mathcal{L}}) \,\, d \mu_{X_{\mathcal{C}} | X_{\mathcal{L}}} 
\,\, \rho_{X_{\mathcal{L}}} (x_{\mathcal{L}})  \,\, d\mu_{X_{\mathcal{L}}}
\end{align}
where $\rho_{X_{\mathcal{C}}|X_{\mathcal{L}}}$ is the conditional pdf on the generalized contour corresponding to $X_{\mathcal{L}}$ while the measure $\mu_{X_{\mathcal{C}} | X_{\mathcal{L}}}$ satisfies
\begin{align} \label{SISC:eq:disint_thm_volume}
\mu_Z(A) = \displaystyle \int_{\pi(A)} \int_{\pi^{-1}(x_{\mathcal{L}})\cap A} d \mu_{X_{\mathcal{C}} | X_{\mathcal{L}}} \,\, d\mu_{X_{\mathcal{L}}}
\end{align}
with $\mu_Z$ and $\mu_{X_{\mathcal{L}}}$ specified.
%
\end{theorem}

\begin{figure}[h!]
		\centering
		\includegraphics[width = 0.4\textwidth]		
						{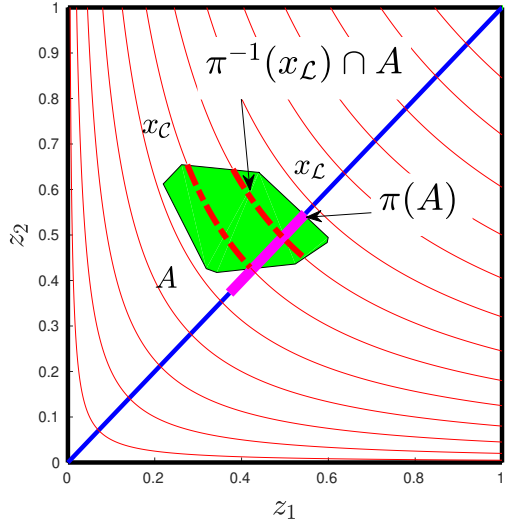}
		\caption{ Illustration of the disintegration theorem. $A$ is represented by the green region while the dashed red curves represent $\pi^{-1}(x_{\mathcal{L}}) \cap A$ and the solid magenta curve represents $\pi(A)$.}
		\label{SISC:fig:DisIntThmEx}
\end{figure}

Figure~\ref{SISC:fig:DisIntThmEx} illustrates the region of integration in \eqref{SISC:eq:disintegratn_thm} for $A\subset \Gamma$.

In Theorem~\ref{SISC:thm:DisintThm}, we have assumed that the disintegration of $P_Z$ into a marginal and conditional family of measures is such that the latter two measures are absolutely continuous with respect to $\mu_{X_{\mathcal{L}}}$ and $\mu_{X_{\mathcal{C}} | X_{\mathcal{L}}}$ to admit the pdfs given by $\rho_{X_{\mathcal{L}}}(x_{\mathcal{L}})$ and $\rho_{X_{\mathcal{C}}|X_{\mathcal{L}}} (x_{\mathcal{C}}|x_{\mathcal{L}})$.  If the probability density functions $\rho_{X_{\mathcal{L}}}$ and $\rho_{X_{\mathcal{C}}|X_{\mathcal{L}}}$ are specified, the pdf of $Z$ can then be estimated. One such computational approach is to partition $\Gamma$ into Voronoi cells and estimate \eqref{SISC:eq:disintegratn_thm} on each partition as in \cite{paper:ButlerGEDW2015}. Due to the bijection between $\mathcal{L}$ and $\mathcal{D}$, the disintegration theorem guarantees that when the pdf of $Z$ constructed in this manner is propagated (pushed forward) through the model, the resulting pdf on $Q$ matches the prescribed $\rho_Q$. This is because $\rho_{X_{\mathcal{L}}}$ is fully specified given knowledge of $\rho_{Q}$ from the relation $\displaystyle \int_K \rho_{X_{\mathcal{L}}} (x_{\mathcal{L}}) \, d \mu_{X_{\mathcal{L}}} = \int_{Q(K)} \rho_Q(q) \, d \mu_Q$ for measurable $K \subset \mathcal{L}$. On the other hand, $\rho_{X_{\mathcal{C}}|X_{\mathcal{L}}}$ cannot be identified solely relying on information from $\rho_Q(q)$, rendering the inverse problem ill-posed. Two distinct choices for $ \rho_{X_{\mathcal{C}}|X_{\mathcal{L}}}$ yield distinct pdfs on $Z$ whose resulting pdf on $Q$ when propagated through the model both match $\rho_Q$. It has been argued that it is reasonable to assume that
$ \rho_{X_{\mathcal{C}}|X_{\mathcal{L}}}$ is uniform over the generalized contour, i.e.

\begin{equation} \label{SISC:eq:ansatz_contour_general}
\rho^{ansatz}_{X_{\mathcal{C}} | X_{\mathcal{L}}} (x_{\mathcal{C}} | x_{\mathcal{L}}) \coloneqq \rho_{X_{\mathcal{C}} | X_{\mathcal{L}}} (x_{\mathcal{C}} | x_{\mathcal{L}}) = \displaystyle  \left( \int_{\pi^{-1} (x_{\mathcal{L}})} \, d \mu_{X_{\mathcal{C}}|X_{\mathcal{ L}}} \right)^{-1} 
\end{equation}
(see \cite[equation 4.5]{paper:ButlerETDW2014}). The ansatz \eqref{SISC:eq:ansatz_contour_general} has been employed with apparently satisfactory results in examples related to recovering the true probability law of $Z$ in PDE models \cite{paper:ButlerETDW2014} as well as in inverse problems that arise in storm-surge applications \cite{paper:ButlerETDW2014}, hydrodynamic models \cite{paper:ButlerGEDW2015}, and groundwater contamination \cite{paper:MattisBDEV2015}. 

To summarize, the method discussed constructs a \emph{consistent} probability measure on $Z$, i.e. the push-forward of this measure matches the observed probability measure on $Q(Z)$. In other words, it constructs a particular pullback measure of an observed measure. It is therefore reasonable to use this method to attempt to ``recover'' the true probability measure/pdf on $Z$ assuming information on $Q(Z)$ only. Of course, this can be successfully accomplished if $m=n$ and if $Q$ is bijective through a standard change of variables. When this is not the case, the ansatz \eqref{SISC:eq:ansatz_contour_general} is likely insufficient if the objective is beyond consistency. 

In order to apply the method above, the measures $\mu_Z, \mu_Q, \mu_{X_{\mathcal{L}}}$ have to be specified. In the next 2 sections, we examine design approaches for these measures as presented in literature \cite{paper:ButlerETDW2014,paper:ButlerGEDW2015}. We also investigate the issues that may arise from application of the ansatz \eqref{SISC:eq:ansatz_contour_general}.

\subsubsection{Pdfs with respect to Lebesgue measures} \label{SISC:subsubsec:DesignLebesgueMeas}

Following Figures 3, 4 and Equations 4.3, 4.4 of \cite{paper:ButlerGEDW2015}, we interpret $\mu_Z, \mu_Q, \mu_{X_{\mathcal{L}}}$ to be chosen as the Lebesgue measure so that $d\mu_Z(z) = dz$, $\mu_Q(q) = dq$, $d\mu_{X_{\mathcal{L}}}(x_{\mathcal{L}}) = dx_{\mathcal{L}}$ and the respective densities are $f_Z, f_Q, f_{X_{\mathcal{L}}}$. Using geometric arguments, it follows from \eqref{SISC:eq:disint_thm_volume} that $\mu_{X_{\mathcal{C}}|X_{\mathcal{ L}}}$ is also Lebesgue so that $d \mu_{X_{\mathcal{C}}|X_{\mathcal{ L}}} = d x_{\mathcal{C}}$. As a consequence, the disintegration theorem \eqref{SISC:eq:disint_thm_volume} can be rewritten as
\begin{align}  \label{SISC:eq:disIntForProbLebesgue}
P_Z(A) = \displaystyle \int_{\pi(A)} \int_{\pi^{-1}(x_{\mathcal{L}})\cap A} f_{X_{\mathcal{C}} | X_{\mathcal{L}}} (x_{\mathcal{C}} | x_{\mathcal{L}}) \,\, d x_{\mathcal{C}}
\,\, f_{X_{\mathcal{L}}} (x_{\mathcal{L}})  \,\, dx_{\mathcal{L}} 
\end{align}
while the ansatz \eqref{SISC:eq:ansatz_contour_general} becomes
\begin{align} \label{SISC:eq:ansatzLebesgue}
f^{ansatz}_{X_{\mathcal{C}} | X_{\mathcal{L}}} (x_{\mathcal{C}} | x_{\mathcal{L}}) = \displaystyle  \left( \int_{\pi^{-1} (x_{\mathcal{L}})} \, d x_{\mathcal{C}} \right)^{-1}.
\end{align}
Note the change in notation to emphasize that the pdfs are now with respect to the Lebesgue measure. 

\textbf{Demonstration of method on an example.} Since all measures in the disintegration theorem have been specified, the method can now be applied to solve stochastic inverse problems. In the following, we consider a simple example to emphasize two points. First, assuming that the pdf on the generalized contours is uniform may not always enable the recovery of the true pdf of $Z$ for any mapping. The calculations carried out that indicate that the true pdf is under/overestimated is supported by Figure~\ref{SISC:fig:arclength_pdf}. Second, the pdf along the generalized contours can have substantial variation for different types of laws on $Z$. Figures~\ref{SISC:fig:contours_and_pdfs} and~\ref{SISC:fig:contourDensGeneral} show how complicated the behavior can be of the pdf on the contour. 

\begin{example} \label{SISC:ex:butler_estep}
Suppose that the inverse problem is independent of the spatial discretization $x$. Consider $Q(Z) = Z_1 \cdot Z_2, \, Z = (Z_1,Z_2)$ where $Z_1,Z_2 \sim U(0,1)$ and are independent. It is shown that given the pdf $f_Q$ of $Q$, (1) the methodology above is unable to recover the probability law of $Z$ and (2) the pdf along the generalized contours of $Q$ are not uniform.
\end{example}
Consider the contours of $Q$ (red thin curves) with the transverse parameterization $\mathcal{L}$ (thick blue line) where $0 \le x_{\mathcal{L}} \le \sqrt{2}$ in Figure~\ref{SISC:fig:new_coord_system}. In order to apply the proposed methodology, suppose that the support $\Gamma \coloneqq Z(\Omega) = [0,1] \times [0,1]$ is known and that $f_Q$ is given. The pdf $f_Q$ can be computed as follows: for $q \in (0,1]$,
\begin{align} \label{SISC:eq:prodPdfDerivation}
P(Q \le q) = \int_{0}^{1} P \left(Z_2 \le \frac{q}{z_1}\right) f_{Z_1} (z_1) \, dz_1 = \int_0^q \, dz_1 + \int_{q}^1 \frac{q}{z_1} \,dz_1 = q - q\log(q)
\end{align}
where $f_{Z_1}$ is the pdf of $Z_1$. Thus, $f_{Q}(q) = -\log(q)$ for $0 < q \le 1$. Given $f_Q$, it remains to determine $f_{X_{\mathcal{L}}}$ and $ f_{X_{\mathcal{C}}|X_{\mathcal{L}}}$ in order to compute \eqref{SISC:eq:disIntForProbLebesgue}. 


\begin{figure}[h!]
		\centering
		\includegraphics[width = 0.5\textwidth]		
						{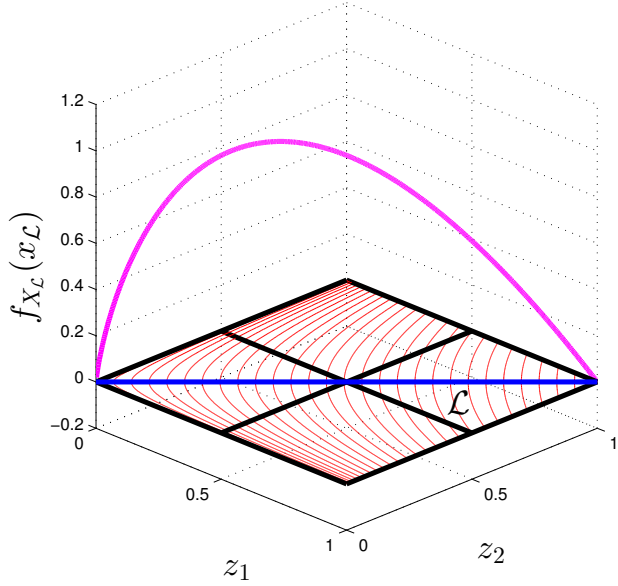}
		\caption{ Plot of the pdf $f_{X_{\mathcal{L}}}$ (thick magenta line) over $\mathcal{L}$  (thick blue line).}
		\label{SISC:fig:new_coord_system_w_pdf}
\end{figure}

As mentioned above, $f_{X_{\mathcal{L}}}$ can be uniquely obtained from $f_Q$. For $x_{\mathcal{L}} \in \mathcal{L}$, the generalized contour $\pi^{-1}(x_{\mathcal{L}})$ corresponding  to $x_{\mathcal{L}}$ passes through the point $(z_1,z_2) = \left(\frac{x_{\mathcal{L}}}{\sqrt{2}},\frac{x_{\mathcal{L}}}{\sqrt{2}}\right)$. Hence, the contour $\pi^{-1}(x_{\mathcal{L}})$ can be parameterized as $z_2 = \frac{x_{\mathcal{L}}^2}{2 z_1}$, $x_{\mathcal{L}}^2/2 \le z_1 \le 1$, and we deduce the relationship $Q = \frac{X_{\mathcal{L}}^2}{2}$.
We therefore have that for $x_{\mathcal{L}} \in (0,\sqrt{2}]$, $P(X_{\mathcal{L}} \le x_{\mathcal{L}}) = P(2Q \le x_{\mathcal{L}}^2) = F_Q(x_{\mathcal{L}}^2 /2)$ which yields
 $$f_{X_{\mathcal{L}}} (x_{\mathcal{L}}) = -x_{\mathcal{L}} \cdot \log(x_{\mathcal{L}}^2/2).$$ 
Figure~\ref{SISC:fig:new_coord_system_w_pdf} shows a plot of $f_{X_{\mathcal{L}}}$ (thick magenta line) over $\mathcal{L}$. Using \eqref{SISC:eq:ansatzLebesgue} and the arc length formula applied to the parameterization of the contour $\pi^{-1}(x_{\mathcal{L}})$ for $x_{\mathcal{L}} \in \mathcal{L}$, we obtain for any $x_{\mathcal{C}} \in \pi^{-1}(x_{\mathcal{L}})$: 

$$f^{ansatz}_{X_{\mathcal{C}}|X_{\mathcal{L}}} (x_{\mathcal{C}}|x_{\mathcal{L}}) = \left( \displaystyle \int_{x_\mathcal{L}^2/2}^1 \sqrt{1 + \frac{x_{\mathcal{L}}^4}{4z_1^4}} \,dz_1 \right)^{-1}.$$

With the components of \eqref{SISC:eq:disIntForProbLebesgue} specified, $P_Z(A)$ can now be approximated for any measurable $A \subset \Gamma = [0,1]^2$. We partition $\Gamma$ into 4 measurable regions of equal area as shown in Figure~\ref{SISC:fig:new_coord_system}, namely the northwest $(A_1)$, northeast $(A_2)$, southwest $(A_3)$, and southeast $(A_4)$ regions. If the proposed methodology is able to recover the true probability law of $Z$ given $f_Q$, we expect the approximations of $P_Z(A_i)$ to be close to 0.25. Numerical calculations reveal that $P_Z(A_i) \approx 0.2886, 0.2459, 0.1770, 0.2886$ for $i = 1,\dots,4$, respectively which shows that $P_Z(A_3)$ in particular underestimates the true probability by a significant amount.

To understand the values obtained for $P_Z(A_i)$, we plot each term in the integrand in \eqref{SISC:eq:disIntForProbLebesgue}, namely $\displaystyle \int_{\pi^{-1}(x_{\mathcal{L}}) \cap A_i} f^{ansatz}_{X_{\mathcal{C}}|X_{\mathcal{L}}} (x_{\mathcal{C}}|x_{\mathcal{L}}) \,dx_{\mathcal{C}}$ for $i=1,\dots,4$ and $f_{X_{\mathcal{L}}}(x_{\mathcal{L}}) $, both as a function of $x_{\mathcal{L}}$ in Figure~\ref{SISC:fig:arclength_pdf}. Each subplot corresponds to each quadrant.  The first term of the integrand amounts to the proportion of the generalized contour inside the quadrant and is displayed with red solid curves. The blue dashed curves meanwhile refer to $f_{x_{\mathcal{L}}}$ for values of $x_{\mathcal{L}}$ contained within the quadrant.

It is worth noting that $P_Z(A_1)$ and $P_Z(A_4)$ are the largest because of the following reasons. Firstly, the regions $A_1$ and $A_4$ are spanned by $0 \le x_{\mathcal{L}} \le 1$ whereas  $f_{X_{\mathcal{L}}}(x_{\mathcal{L}}) $ attains its maximum in the neighborhood $0.2 \le x_{\mathcal{L}} \le 0.8$. Secondly, these two quadrants possess more contours than regions $A_2$ and $A_3$ because every contour passing through the 2 former regions always passes through one of the latter regions. This implies that $A_1$ and $A_4$ comprise a wider range of values of $x_{\mathcal{L}}$ as evidenced by the domains of the corresponding subplots in Figure~\ref{SISC:fig:arclength_pdf}.

\begin{figure}[h!]
		\centering
		\includegraphics[width = 0.8\textwidth]		
						{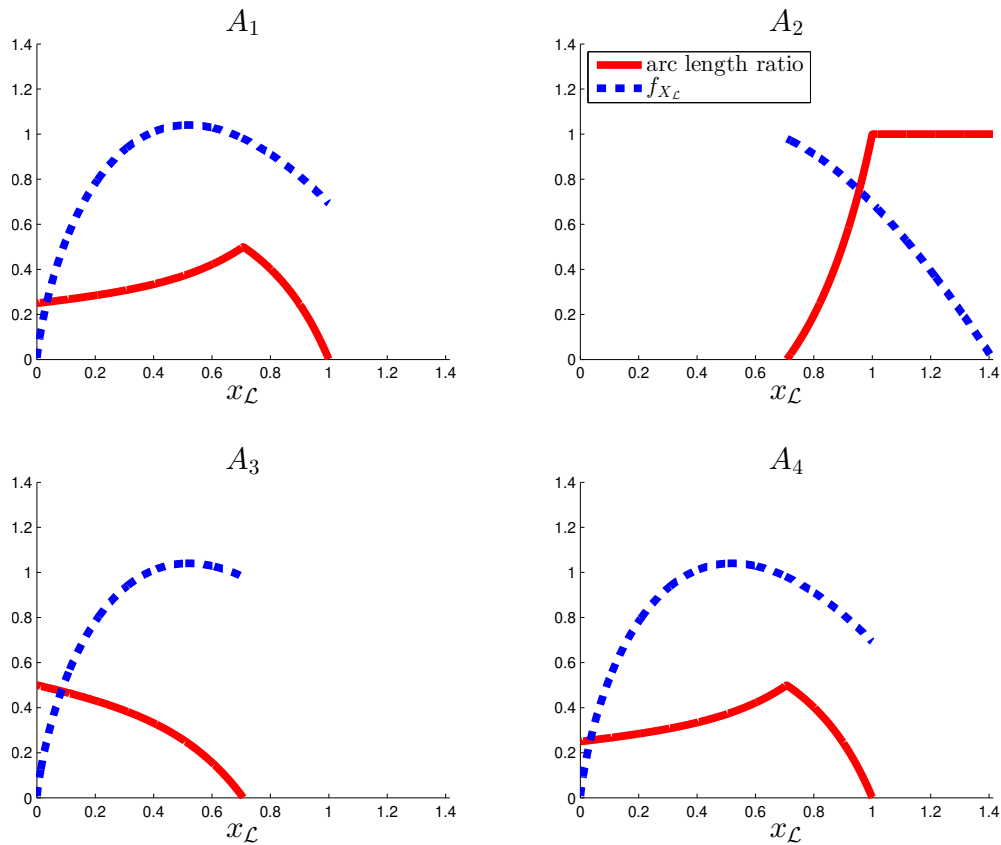}
		\caption{Plot of the integrand terms in \eqref{SISC:eq:disIntForProbLebesgue} as a function of  $x_{\mathcal{L}}$ for each $A_i$. The blue dashed curve is $f_{X_{\mathcal{L}}}$ while the red solid curve is the proportion of each contour contained in $A_i$.}
		\label{SISC:fig:arclength_pdf}
\end{figure}

Finally, we numerically demonstrate that for this example, the actual conditional pdf along the generalized contours is not necessarily uniform (see Appendix~\ref{SISC:appendix:pdfContour} for the methodology). The actual conditional pdf results from disintegrating the measure $P_Z$ associated with the uniform distribution on $Z$. The left panel of Figure~\ref{SISC:fig:contours_and_pdfs} shows three contours of $Q$, the middle panel shows the corresponding actual pdfs along each contour, whereas for comparison, the right panel shows the pdf according to the ansatz \eqref{SISC:eq:ansatzLebesgue}. As the concavity of the contours increases, the actual pdf along the contour becomes less uniform. We see that the pdf on each contour resulting from the disintegration is not necessarily uniform even though $Z$ is uniformly distributed on $\Gamma$. In addition, Figure~\ref{SISC:fig:contours_and_pdfs} corroborates the results in Figure~\ref{SISC:fig:arclength_pdf}: the majority of the contours residing in $A_3$ exhibit high concavity which lead to the underestimation of $P_Z(A_3)$.  If instead the $f_{X_{\mathcal{C}} \,| \,X_{\mathcal{L}}} (x_{\mathcal{C}}|x_{\mathcal{L}})$ as obtained in the middle panel of Figure~\ref{SISC:fig:contours_and_pdfs} were used in computing $P_Z(A)$ in \eqref{SISC:eq:disIntForProbLebesgue}, the resulting $Z$ would be uniformly distributed. Example~\ref{SISC:ex:butler_estep} therefore underscores the point made at the beginning of this section that {\em if the objective is to recover the true pdf on $Z$, then the ansatz \eqref{SISC:eq:ansatzLebesgue}, which only ensures that the constructed measure/density is a pullback measure, may be insufficient to recover the structure of the true pdf in directions not informed by the data}.

\begin{figure}[h!]
		\centering
		\includegraphics[width = 1\textwidth]		
						{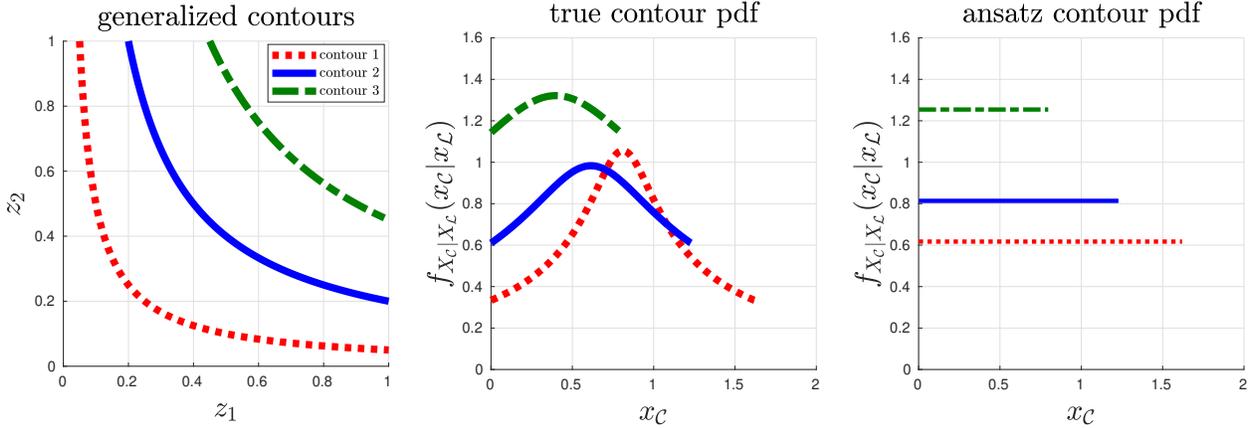}
		\caption{Left panel: selected contours of $Q$ in Example~\ref{SISC:ex:butler_estep}. Middle panel: corresponding actual pdf along the contour. Right panel: corresponding pdf using the ansatz~\eqref{SISC:eq:ansatzLebesgue}.}
		\label{SISC:fig:contours_and_pdfs}
\end{figure}

Furthermore, we evaluate how \eqref{SISC:eq:ansatzLebesgue} fares as a means of regularizing across other pdfs on $X_{\mathcal{C}}| X_{\mathcal{L}}$ in the absence of information on $Z$. Suppose instead that $Z_1,Z_2$ are independent with $Z_1 \sim Beta(\nu_1,\nu_2)$ and $Z_2 \sim Beta(\tau_1,\tau_2)$. Notice that the pdf of $Z$ in Example~\ref{SISC:ex:butler_estep} is a special case with $\nu_1,\nu_2,\tau_1,\tau_2$ all being equal to 1. The 3 subplots of Figure~\ref{SISC:fig:contourDensGeneral} show the conditional pdf along the solid blue contour in Figure~\ref{SISC:fig:contours_and_pdfs} (contour 2) for different combinations of the parameters for $Z_1,Z_2$. The plots reveal that the pdfs on the contours can be very complex and suggest that using the pdf along the contour as the only means of regularizing against other solutions to this inverse problem may be insufficient.

\begin{figure}[h!]
		\centering
		\includegraphics[width = 1\textwidth]		
						{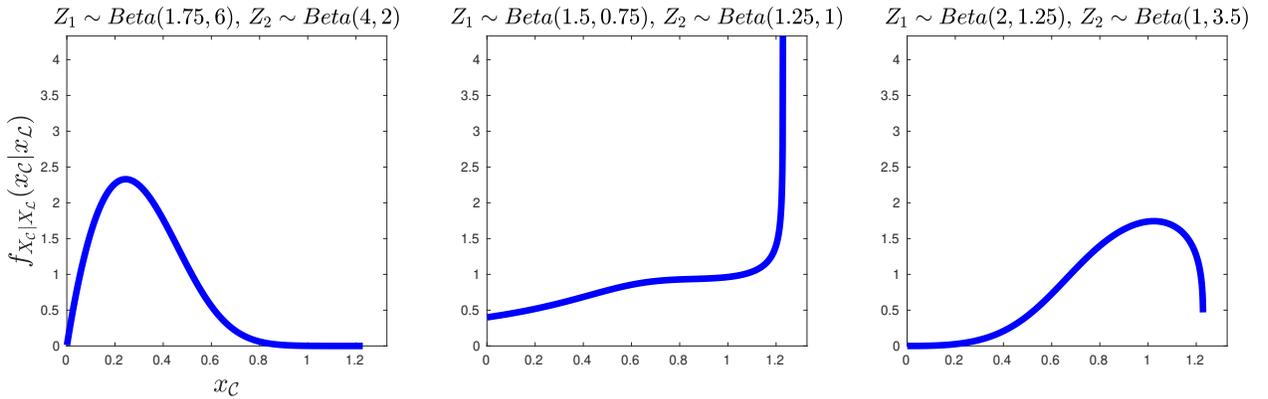}
		\caption{Conditional pdf $X_{\mathcal{C}}|X_{\mathcal{L}}$ on contour 2 in Figure~\ref{SISC:fig:contours_and_pdfs} where $Z_1 \sim Beta(\nu_1,\nu_2), Z_2 \sim Beta(\tau_1,\tau_2)$. }
		\label{SISC:fig:contourDensGeneral}
\end{figure}

\subsubsection{Pdfs with respect to non-Lebesgue measures}\label{SISC:subsubsec:DesignNonLebesgue}

In contrast to Section~\ref{SISC:subsubsec:DesignLebesgueMeas}, \cite{paper:ButlerETDW2014} chose the measures in \eqref{SISC:eq:disintegratn_thm} as follows: $\mu_Z$ is taken to be Lebesgue, $\mu_Q$ is defined to be the pushforward measure of $\mu_Z$ through $Q$ (i.e. for measurable $B \subset \mathcal{D}$, $\mu_Q (B) = \mu_Z(Q^{-1}(B))$), while $\mu_{X_{\mathcal{L}}}$ is computed as
$\mu_{X_{\mathcal{L}}} (K) = \mu_Q(Q(A))$ where $K = \pi(A)$. Note that $\mu_Q$ may not be Lebesgue especially if $Q$ is nonlinear which implies that the same holds for $\mu_{X_{\mathcal{L}}} $. As before, $\mu_{X_{\mathcal{C}} | X_{\mathcal{L}}}$ results from \eqref{SISC:eq:disint_thm_volume} since $\mu_Z,\mu_{X_{\mathcal{L}}}$ are now specified; however, if  $\mu_{X_{\mathcal{L}}}$ is not be the Lebesgue measure, then neither is $\mu_{X_{\mathcal{C}} | X_{\mathcal{L}}}$.

To further distinguish this design approach from that of Section~\ref{SISC:subsubsec:DesignLebesgueMeas}, we emphasize that for fixed $x_{\mathcal{L}}$, $\rho^{ansatz}_{X_{\mathcal{C}} | X_{\mathcal{L}}} (x_{\mathcal{C}} | x_{\mathcal{L}})$ is constant with respect to $\mu_{X_{\mathcal{C}} | X_{\mathcal{L}}}$. If we were to express this pdf with respect to the Lebesgue measure, i.e. we seek $f^{ansatz}_{X_{\mathcal{C}} | X_{\mathcal{L}}} (x_{\mathcal{C}} | x_{\mathcal{L}})$ such that $\rho^{ansatz}_{X_{\mathcal{C}} | X_{\mathcal{L}}} (x_{\mathcal{C}} | x_{\mathcal{L}}) \, d \mu_{X_{\mathcal{C}} | X_{\mathcal{L}}} = f^{ansatz}_{X_{\mathcal{C}} | X_{\mathcal{L}}} (x_{\mathcal{C}} | x_{\mathcal{L}}) \, d x_{\mathcal{C}}$, then $f^{ansatz}_{X_{\mathcal{C}} | X_{\mathcal{L}}} (x_{\mathcal{C}} | x_{\mathcal{L}})$ may not be constant unlike in \eqref{SISC:eq:ansatzLebesgue}. In addition, this approach offers a computationally efficient approximation for the case when $Z,Q(Z)$ are high-dimensional that does not require explicitly constructing $\mu_{X_{\mathcal{C}} | X_{\mathcal{L}}}$ from \eqref{SISC:eq:disint_thm_volume} \cite[Algorithm 1]{paper:ButlerETDW2014}. The publicly available code\footnote{https://github.com/UT-CHG/BET} is based on this design.

\textbf{Demonstration of method on an example.} Even with this selection of the measures in \eqref{SISC:eq:disintegratn_thm}, it is shown that the ansatz \eqref{SISC:eq:ansatz_contour_general} may still be unable to recover the true law of $Z$. We show this in the next example where the true distribution of $Z$ in Example~\ref{SISC:ex:butler_estep} is modified. The established results in 
Figures~\ref{SISC:fig:contours_and_pdfs} and~\ref{SISC:fig:contourDensGeneral} of Section~\ref{SISC:subsubsec:DesignLebesgueMeas} will be invoked to aid in the discussion.

\begin{example} \label{SISC:ex:butler_estep_NonLebesgue}
We revisit the setup in Example~\ref{SISC:ex:butler_estep} but with the true distribution of $Z$ altered to $Z_1 \sim Beta(\nu_1,\nu_2), Z_2 \sim Beta(\tau_1,\tau_2)$, independent. Given the pdf $\rho_Q$ of $Q$, it is shown that the probability law of $Z$ cannot be recovered with the choice of ansatz \eqref{SISC:eq:ansatz_contour_general}.
\end{example}

In order to apply the disintegration theorem \eqref{SISC:eq:disintegratn_thm}, $\mu_{X_{\mathcal{C}} | X_{\mathcal{L}}}$ needs to be first identified from \eqref{SISC:eq:disint_thm_volume}. Let $P^{Unif}_Z$ be the probability measure of the uniform distribution on $Z$ where $Z_1,Z_2$ are independent. For measurable $A \subset \Gamma$, \eqref{SISC:eq:disIntForProbLebesgue} and the results in Section~\ref{SISC:subsubsec:DesignLebesgueMeas} yield the disintegration
\begin{align} \label{SISC:eq:DisIntUnifDist}
P_Z^{Unif}(A) = \displaystyle \int_{\pi(A)} \int_{\pi^{-1}(x_{\mathcal{L}})\cap A} f^{Unif}_{X_{\mathcal{C}} | X_{\mathcal{L}}} (x_{\mathcal{C}} | x_{\mathcal{L}}) \,\, f^{Unif}_{X_{\mathcal{L}}} (x_{\mathcal{L}}) \,\, d x_{\mathcal{C}}
\,\,  dx_{\mathcal{L}}
\end{align}
where $f^{Unif}_{X_{\mathcal{C}} | X_{\mathcal{L}}}$ are the pdfs displayed in the middle panel of Figure~\ref{SISC:fig:contours_and_pdfs}. Since $\mu_Z$ is Lebesgue, $P_Z^{Unif}(A) = \mu_Z(A)$ which implies that
\begin{align} \label{SISC:eq:VolumeMeasureOnContour}
d \mu_{X_{\mathcal{C}} | X_{\mathcal{L}}} = f^{Unif}_{X_{\mathcal{C}} | X_{\mathcal{L}}} (x_{\mathcal{C}} | x_{\mathcal{L}}) \,\, d x_{\mathcal{C}}
\end{align}
by comparing \eqref{SISC:eq:DisIntUnifDist} with \eqref{SISC:eq:disint_thm_volume}.

We now return to the setup in Example~\ref{SISC:ex:butler_estep_NonLebesgue} where $Z_1,Z_2$ are actually independent beta random variables. Denote by $P_Z^{Beta}$ the probability measure on $Z$ corresponding to its true law. From the results in Section~\ref{SISC:subsubsec:DesignLebesgueMeas}, its actual disintegration is
\begin{align} \label{SISC:eq:DisIntBetaDist}
P_Z^{Beta}(A) = \displaystyle \int_{\pi(A)} \int_{\pi^{-1}(x_{\mathcal{L}})\cap A} f^{Beta}_{X_{\mathcal{C}} | X_{\mathcal{L}}} (x_{\mathcal{C}} | x_{\mathcal{L}}) \,\, f^{Beta}_{X_{\mathcal{L}}} (x_{\mathcal{L}}) \,\, d x_{\mathcal{C}}
\,\,  dx_{\mathcal{L}}
\end{align}
where $f^{Beta}_{X_{\mathcal{C}} | X_{\mathcal{L}}}$ is obtained similarly as the pdfs plotted in the panels of Figure~\ref{SISC:fig:contourDensGeneral}. On the other hand, applying the ansatz \eqref{SISC:eq:ansatz_contour_general} to \eqref{SISC:eq:disintegratn_thm},  the inverse problem solution is 
\begin{align} \label{SISC:eq:ansatzSolnBeta}
P_Z^{ansatz}(A) =  \displaystyle \int_{\pi(A)} \int_{\pi^{-1}(x_{\mathcal{L}})\cap A} \rho^{ansatz}_{X_{\mathcal{C}} | X_{\mathcal{L}}} (x_{\mathcal{C}} | x_{\mathcal{L}}) \,\, \rho_{X_{\mathcal{L}}} (x_{\mathcal{L}}) \,\, d \mu_{X_{\mathcal{C}} | X_{\mathcal{L}}}   \,\, d\mu_{X_{\mathcal{L}}}
\end{align}
where $\rho^{ansatz}_{X_{\mathcal{C}} | X_{\mathcal{L}}} (x_{\mathcal{C}} | x_{\mathcal{L}})$ is a constant for fixed $x_{\mathcal{L}}$. If \eqref{SISC:eq:ansatzSolnBeta} is able to recover the true law of $Z$, it must be that $P^{Beta}_Z(A)= P^{ansatz}_Z(A)$ for any measurable $A \subset \Gamma$. Since the probability measure on $X_{\mathcal{L}}$ is the same regardless of the choice of $\mu_{X_{\mathcal{L}}}$, it follows that $\rho_{X_{\mathcal{L}}} (x_{\mathcal{L}}) \,\, d\mu_{X_{\mathcal{L}}} = f^{Beta}_{X_{\mathcal{L}}} (x_{\mathcal{L}}) \,\, dx_{\mathcal{L}}.$ Hence, by comparing \eqref{SISC:eq:ansatzSolnBeta} with \eqref{SISC:eq:DisIntBetaDist}, the ansatz is able to recover the true law if and only if 
\begin{align} \label{SISC:eq:conditionForRecovery}
\rho^{ansatz}_{X_{\mathcal{C}} | X_{\mathcal{L}}} (x_{\mathcal{C}} | x_{\mathcal{L}}) \,\, d \mu_{X_{\mathcal{C}} | X_{\mathcal{L}}} = f^{Beta}_{X_{\mathcal{C}} | X_{\mathcal{L}}} (x_{\mathcal{C}} | x_{\mathcal{L}}) \,\, d x_{\mathcal{C}}  
\end{align}
which is equivalent to 
\begin{align} \label{SISC:eq:IFFcondition}
\rho^{ansatz}_{X_{\mathcal{C}} | X_{\mathcal{L}}} (x_{\mathcal{C}} | x_{\mathcal{L}}) = \frac{f^{Beta}_{X_{\mathcal{C}} | X_{\mathcal{L}}} (x_{\mathcal{C}} | x_{\mathcal{L}})}{f^{Unif}_{X_{\mathcal{C}} | X_{\mathcal{L}}} (x_{\mathcal{C}} | x_{\mathcal{L}})}
\end{align}
using \eqref{SISC:eq:VolumeMeasureOnContour}. 
%

We now proceed by contradiction. Consider contour 2 (solid blue) in the left panel of Figure~\ref{SISC:fig:contours_and_pdfs} and let $\nu_1,\nu_2,\tau_1,\tau_2$ be any of the values utilized in the plots of Figure~\ref{SISC:fig:contourDensGeneral}. The pdf $f^{Beta}_{X_{\mathcal{C}} | X_{\mathcal{L}}} (x_{\mathcal{C}} | x_{\mathcal{L}})$ on this contour could then be any of these plotted pdfs. It then follows from the middle panel of Figure~\ref{SISC:fig:contours_and_pdfs} and each of the panels of Figure~\ref{SISC:fig:contourDensGeneral} that $\frac{f^{Beta}_{X_{\mathcal{C}} | X_{\mathcal{L}}} (x_{\mathcal{C}} | x_{\mathcal{L}})}{f^{Unif}_{X_{\mathcal{C}} | X_{\mathcal{L}}} (x_{\mathcal{C}} | x_{\mathcal{L}})}$ is not constant for $x_{\mathcal{L}}$ fixed which contradicts the assumption on $\rho^{ansatz}_{X_{\mathcal{C}} | X_{\mathcal{L}}} (x_{\mathcal{C}} | x_{\mathcal{L}})$. The ansatz is therefore unable to recover the true law of $Z$.

Despite the examples presented in Sections~\ref{SISC:subsubsec:DesignLebesgueMeas} and~\ref{SISC:subsubsec:DesignNonLebesgue}, the methodology in \cite{paper:BreidtBE2011, paper:ButlerETDW2014, paper:ButlerGEDW2015} can still be useful in physical applications in that it can serve as a first model for the unknown pdf of $Z$. This can be later tuned or scrutinized for plausibility depending on available information on $Z$. 

\subsubsection{Incorporating prior information} \label{SISC:sec:ConsistentBayesian}

An alternative to the method in Section~\ref{SISC:subsubsec:BETreview} under the same specifications on the inverse problem has been developed in \cite{paper:ButlerJW2018}. It assumes that some information on $Z$ is available in the form of a prior pdf $\rho_Z^{prior}(Z)$. Analogous to Bayes' theorem, the solution to the inverse problem is a posterior pdf $\rho_Z^{post}$ on $Z$ that is constructed as

$$\rho_Z^{post}(Z) = \rho_Z^{prior}(Z) \cdot \frac{\rho_Q(Q(Z))}{\rho_Q^{Q(prior)}(Q(Z))}$$
for $Z \in \Gamma$ where $\rho_Q(\cdot)$ is the given pdf of $Q$ while $\rho_Q^{Q(prior)}(\cdot)$ is the pdf of $Q$ that is obtained by propagating $\rho_Z^{prior}$ through the model. This solution was derived using the disintegration theorem based on conditional densities which is more general than Theorem~\ref{SISC:thm:DisintThm} based on generalized contours. Although this method does not explicitly deal with contours, it is related to the method in Section~\ref{SISC:subsec:BETMethod} in that $\rho_Z^{prior}(Z)$ implies a pdf $\rho_{X_{\mathcal{C}}|X_{\mathcal{L}}} (x_{\mathcal{C}}|x_{\mathcal{L}})$ on the generalized contours that is not necessarily uniform. We note that a sufficient condition for this method to recover the true pdf $\rho_Z$ on $Z$ would be if $\rho_Z^{prior} = \rho_Z$, signifying no gain in information. A more general sufficient condition only requires that the conditional pdf along the contours arising from the disintegration of $\rho_Z^{prior}$ and the true pdf $\rho_Z$ need to be equal \cite[Sections 3, 7.3]{paper:ButlerJW2018}. The methodology proposed in \cite{paper:ButlerJW2018} only stresses the need for additional information to be specified on $Z$ in order to solve the inverse problem.


\subsection{Parametric representations of the unknown random field}\label{SISC:subsec:VanWyk}

This section elaborates on the second approach \cite{paper:BorggaardV2015} which estimates the coefficient field of a differential equation given observations of the solution field. Section~\ref{SISC:subsubsec:VanWykReview} reviews the methodology and establishes its similarity with the above inverse problem. Section~\ref{SISC:subsubsec:VanWykExample} meanwhile explores two examples utilizing this method. It is shown that the probability law and the truncation level of the random variables arising from the Karhunen-Lo\`eve expansion of the solution field may be inadequate to characterize the coefficient field.

\subsubsection{Review of methodology}\label{SISC:subsubsec:VanWykReview}

Consider the stochastic equation $\mathcal{L}(U(x,\omega)) = 0$ defined on the probability space $(\Omega,\mathcal{F},P)$ for $\omega \in \Omega, x \in D$ where the operator $\mathcal{L}$ characterizes a stochastic differential equation that depends on the random field $A(x,\omega)$. The inverse problem tackled by \cite{paper:NarayananZ2004,paper:ZabarasG2008,paper:BorggaardV2015,paper:DesceliersGS2006} arising from this set-up is:

\begin{quote}
	Given observations $\hat{U}$ of the solution $U(x,\omega)$, estimate the unknown field $A(x,\omega)$.
\end{quote}

Although this problem appears different from the one posed in Section~\ref{SISC:sec:Intro}, they are in fact conceptually identical. In practice, the spatial domain is discretized so that the random fields $A$ and $U$ are represented as random vectors characterized by the finite-dimensional distributions of the random fields. The inverse problem is now reminiscent of the problem above. The use of the finite-dimensional distribution of $A(x,\omega)$ to characterize the field itself can be justified under the mild assumption that $A(x,\omega)$ has almost surely continuous sample paths or that it satisfies the H\"older continuity condition \cite[Theorem 3.1]{book:Grigoriu2012}; see \cite{book:Grigoriu2012,book:WongH1985} for more details. In what follows, it will be assumed that either assumption on $A$ holds.


Despite this connection, proposed methodologies \cite{paper:NarayananZ2004,paper:ZabarasG2008,paper:BorggaardV2015,paper:DesceliersGS2006} estimate the unknown random field by finding a finite-dimensional noise approximation, that is, $A(x,\omega) \approx A(x,Z(\omega))$ where $Z$ is a random vector. Essentially, these methodologies construct a parametric representation of the unknown field $A$ in which the law, the dimension of $Z$, and the functional form of $A(x,Z)$ are to be determined. In \cite{paper:DesceliersGS2006}, each observed sample of $U$ is used to acquire samples of $A$ via optimization; these samples of $A$ then serve to calibrate a polynomial chaos expansion (PCE) \cite{book:Grigoriu2012} for $A$. This optimization procedure may result in non-unique global minima yet its implications on the fitted PCE were not addressed. In \cite{paper:NarayananZ2004}, observed samples of $U$ are utilized to obtain a truncated PCE of $U$. The unknown $A$ is then expressed in terms of this stochastic basis. The limitations of PCE are well known \cite{paper:FieldG2007,paper:FieldG2005} and furthermore, as the stochastic basis for $U$ is used as the stochastic basis for the unknown $A$, the truncation level of the PCE for $U$ can be insufficient for the PCE of $A$. The following illustrations tackles some of these issues.

In order to address the limitations of a PCE model for $A$, \cite{paper:ZabarasG2008,paper:BorggaardV2015} proposed to express the unknown random field as a sparse grid representation $\widetilde{A}^{N}(x,Y)$ of level $N$. This is represented as
\begin{align} \label{SISC:eq:SPGrid}
	\widetilde{A}^N(x,y) = \sum_{j=1}^N v(x,\mathbf{y}_j) \psi_j(y), \,\,\, x \in \mathbb{R}^d, \,\,\, y \in \Gamma \subset \mathbb{R}^k
\end{align}
for $\Gamma$ bounded. Here, $Y \in \Gamma$ is a random vector whose dimension is not necessarily identical to that of $Z$ and whose law has to be specified. In addition, $\{\mathbf{y}_j\}_{j=1}^N \subset \Gamma$ are the $N$ sparse grid nodes, $\psi_j(y)$ are specified interpolating functions on $\Gamma$, while $v(x,\mathbf{y}_j)$ is a deterministic function of $x$. It only remains to address the following: How must the dimension and the law of the random vector $Y$ be specified for the nodal values $v(x,\mathbf{y}_j)$ to be approximated through optimization?

Both \cite{paper:ZabarasG2008,paper:BorggaardV2015} are similar in that they employ the finite-dimensional model \eqref{SISC:eq:SPGrid} for $A$ yet they differ in how the dimension and the law of $Y$ are prescribed. In \cite{paper:ZabarasG2008}, the dimension of $Y$ is postulated to be $k=1,2$ in the numerical examples while the choice of law for $Y$ is downplayed. This was demonstrated with a numerical example in which various distributions, with various support, were chosen for $Y \in \mathbb{R}^1$ to compute moments of $\widetilde{A}^N(x,Y)$. In contrast, \cite{paper:BorggaardV2015} pursued an approach in selecting the dimension and law of $Y$ based on the Karhunen-Lo\`eve (KL) expansion \cite{book:Grigoriu2012} and is detailed as follows. Consider the elliptic system 
\begin{align} \label{SISC:eq:elliptic_eqn}
-\nabla \cdot (A(x,\omega) \nabla U(x,\omega)) = f(x), \,\,\,\,\, U(x,\omega) = 0 \,\,\, \text{on} \,\,\, \partial D
\end{align}
where $x\in D \subset \mathbb{R}^d, \,\, \omega \in \Omega$. Given observations $\hat{U}$ (possibly noisy) of $U$, the objective is
to approximate the unknown field $A$ by solving the optimization problem
\begin{align} \label{SISC:eq:objFun}
\text{min} \,\, \frac{1}{2} \|U - \hat{U}\|^2 + \frac{\beta}{2} \|A\|^2
\end{align}
for $U$ and $A$ under the constraint that they satisfy \eqref{SISC:eq:elliptic_eqn} and that $A$ satisfies ellipticity constraints. The second term in \eqref{SISC:eq:objFun} regularizes the solution $A$ through the parameter $\beta > 0$ while the norms in \eqref{SISC:eq:objFun} are formed using a tensor product of norms for Sobolev spaces and the $L^2(\Omega)$ norm for the probability space. A KL expansion of $\hat{U}$ is then performed to obtain for $\omega \in \Omega$
\begin{align} \label{SISC:eq:KLsamples}
\hat{U}(x,\omega) = \hat{u}_0(x) + \sum_{k=1}^{\infty} \sqrt{\hat{\lambda}_k} \hat{\phi}_k(x) Y_k(\omega)
\end{align}
where $E[Y_k] = 0, \, E[Y_k^2] = 1 \,\, \forall k$, $E[Y_k Y_j] = 0$ for $k\neq j$ and $\hat{\lambda}_k,\, \hat{\phi}_k(x)$ are the eigenvalues and eigenfunctions of the covariance function of $\hat{U}(x,\omega)$. It is assumed \cite[p. 10]{paper:BorggaardV2015} that $\{Y_i\}_{i=1}^{\infty}$ form a basis for $L^2(\Omega)$ in the sense that every random variable with finite variance can be expressed as a linear combination of $\{Y_i\}_{i=1}^{\infty}$. As a consequence, since the optimal random field $A^*$ for $A$ in \eqref{SISC:eq:objFun} satisfies $A^*(x,\cdot) \in L^2(\Omega) \,\, \forall x \in D$,  $A^*$ can be written as
\begin{align} \label{SISC:eq:optimalRF}
A^*(x,\omega) = a_0(x) + \sum_{k=1}^{\infty} a_{k}(x) Y_k(\omega). 
\end{align} 
The deterministic functions $\{a_k(x)\}_{k=0}^{\infty}$ are determined by solving \eqref{SISC:eq:objFun} while $\{Y_k(\omega)\}_{k=1}^{\infty}$ are obtained from \eqref{SISC:eq:KLsamples}. Due to the dependence of $A^*$ on an infinite number of random variables, \eqref{SISC:eq:optimalRF} is usually referred to as the solution to the infinite-dimensional problem.

For practical numerical implementation, the KL expansion in \eqref{SISC:eq:KLsamples} is truncated to only consider $\{Y_k\}_{k=1}^M, \, M \ll \infty $ based on the decay of $\lambda_k$. The optimal solution then takes on the form $A^\dagger(x,\omega) = A^\dagger(x,Y_1(\omega),\dots,Y_M(\omega))$, i.e. the optimal $A^\dagger$ that is sought from \eqref{SISC:eq:objFun} is a function of $Y_1,\dots,Y_M$ only. $A^\dagger$ is typically referred to as the solution to the finite noise problem. It was shown in \cite{paper:BorggaardV2015} that the sequence of minimizers of the finite noise problem has a subsequence that converges weakly to a minimizer of the infinite-dimensional problem under the assumption \eqref{SISC:eq:optimalRF}. Due to this truncation, $A^\dagger(x,Y_1(\omega),\dots,Y_M(\omega))$ is not expected to be linear in $\{Y_k\}_{k=1}^M$ as in \eqref{SISC:eq:optimalRF}; as such, a sparse grid representation \eqref{SISC:eq:SPGrid} for $A^\dagger(x,Y_1(\omega),\dots,Y_M(\omega))$ using linear hat functions for $\psi_j$  was considered in \cite{paper:BorggaardV2015} to accommodate smoothness conditions on $A^\dagger$ as a function of $Y_1,\dots,Y_M$. 

In summary, \cite{paper:BorggaardV2015} parameterized the unknown field $A(x,\omega)$ by a random vector $Y$ whose dimension and distribution are based on the  random variables arising from the truncated KL expansion of the observed samples $\hat{U}$. The approach pursued by \cite{paper:BorggaardV2015} in specifying the random vector $Y$ is reasonable because it is shown in \cite{paper:BabuskaTZ2005} that if $A$ in \eqref{SISC:eq:elliptic_eqn} depends on $Y$, the response $U$ is analytical in $Y$. But challenges may surface from this approach which we demonstrate through examples in the next section. We clarify that the issue does not lie with the use of a sparse grid approximation for $A$ but with the choice of the random vector $Y$ used to construct $A^*(x,\omega)$ in \eqref{SISC:eq:optimalRF}.

\subsubsection{Demonstration of method on examples}\label{SISC:subsubsec:VanWykExample}

This section investigates the implications of the approach in \cite{paper:BorggaardV2015}. First, unlike the PCE, the infinite set of random variables $\{Y_k(\omega)\}_{k=1}^{\infty}$ in \eqref{SISC:eq:KLsamples} do not form a stochastic basis for $L^2(\Omega)$. For this reason, if $A \in L^2(\Omega)$, the truncated set $\{Y_k(\omega)\}_{k=1}^M$ might not even be adequate to characterize the true field $A$ because it is unlikely that $A \in L^2(\Omega) \cap \text{span}(\{Y_k\}_{k=1}^M)$, as Example~\ref{SISC:ex:vanWyk1} will clarify.  Second, even if it can be analytically shown that the random variables $\{Y_k\}_{k=1}^{\infty}$ in the KL expansion of $U$ also characterize $A$, the truncation level employed for the practical implementation of $U$ might not be sufficient for the implementation of $A$. Example~\ref{SISC:ex:vanWyk2} lends support to this claim.

\begin{example} \label{SISC:ex:vanWyk1}
Consider the stochastic ODE for $x \in [0,1], \, \omega \in \Omega$:
\begin{align} \label{SISC:eq:ODEexample}
	-\frac{d}{dx}(A(x,\omega) \cdot \frac{d}{dx}U(x,\omega)) = 0, \,\,\,\,\,\, U(0,\omega) = 0, \,\,\,\,\,\, U(1,\omega) = \int_0^1 \frac{1}{A(y,\omega)} \,dy. 	
\end{align}
The true coefficient random field is modeled as a translation process \cite{book:Grigoriu2012} $A(x,\omega) = \alpha + (\beta - \alpha) \cdot F^{-1}_{beta}(\Phi(G(x,\omega)))$ where $\alpha = 4, \beta = 20$, $F^{-1}_{beta}$ is the inverse cumulative distribution function (cdf) of $Beta(1,3)$, $\Phi$ is the cdf of $N(0,1)$, and $G(x,\omega)$ is a zero mean, unit variance stationary Gaussian process with Mat\'ern covariance, i.e. $E[G(s,\cdot) G(t,\cdot)] = \frac{2^{1-\nu}}{\Gamma(\nu)} \left(\frac{\sqrt{2\nu} |s-t|}{\ell}\right)^{\nu} K_{\nu}\left(\frac{\sqrt{2\nu} |s-t|}{\ell}\right), \,\, s,t \in [0,1]$ where $K_{\nu}$ is the modified Bessel function of the second kind and $\nu = \frac{5}{2}, \ell = 0.03$. To formulate the inverse problem, we approximate $A(x,\omega)$ using observed samples of $U(x,\omega)$. It is shown that the random variables arising from the KL expansion of $U(x,\omega)$ are inadequate to characterize $A(x,\omega)$.  
\end{example}

\begin{figure}[h!]
		\centering
		\includegraphics[width = 0.8\textwidth]		
						{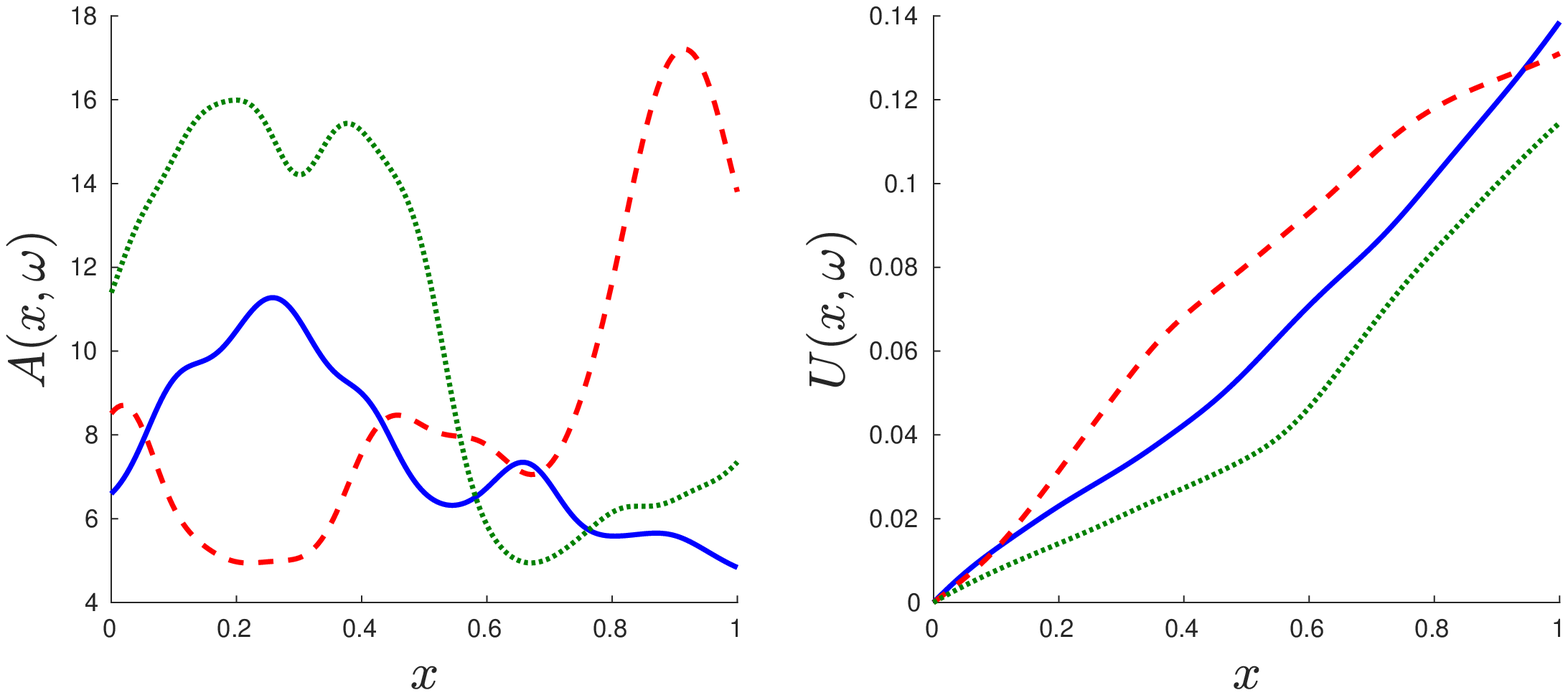}
		\caption{Left panel: samples of $A(x,\omega)$. Right panel: corresponding samples of $U(x,\omega)$ via \eqref{SISC:eq:ODEexample}.}
		\label{SISC:fig:SamplesAndResponse}
\end{figure}

\begin{figure}[h!]
		\centering
		\includegraphics[width = 0.9\textwidth]		
						{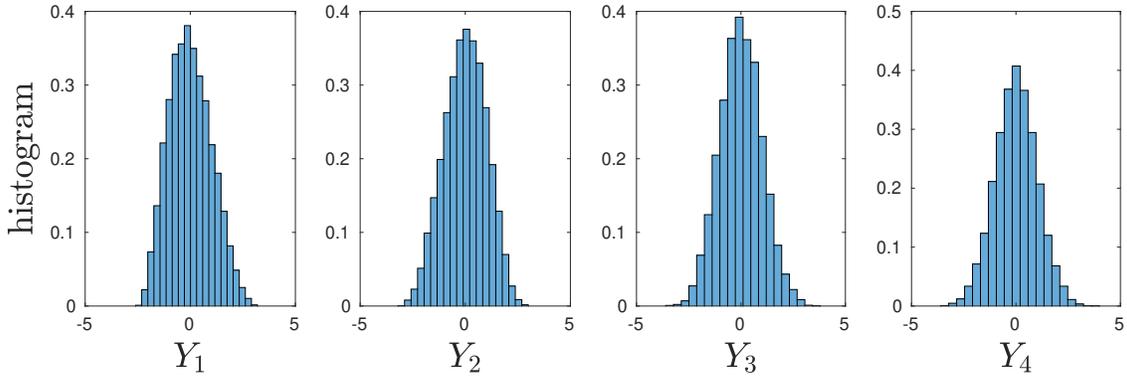}
		\caption{Histograms of the first 4 random variables in the KL expansion of $U$.}
		\label{SISC:fig:KLRVhist}
\end{figure}

We generate 10000 noiseless samples of $U(x,\omega)$ by solving \eqref{SISC:eq:ODEexample} for each sample of $A(x,\omega)$. The analytical solution is given by $U(x,\omega) = \displaystyle \int_0^x \frac{1}{A(y,\omega)} \, dy = \displaystyle \int_0^x B(y,\omega) \, dy$ where $B(y,\omega) \coloneqq \frac{1}{A(y,\omega)}$. A KL expansion of $U(x,\omega)$ is then performed to obtain \eqref{SISC:eq:KLsamples}. Figure~\ref{SISC:fig:SamplesAndResponse} shows samples of $A(x,\omega)$ together with their respective samples $U(x,\omega)$ while Figure~\ref{SISC:fig:KLRVhist} displays histograms of the samples $Y_k$ in the KL expansion of $U(x,\omega)$ corresponding to the four largest eigenvalues. Note that the integral equations that need to be solved to obtain the eigenvalues and eigenfunctions of the covariance function of $U(x,\omega)$ have to be discretized solely for numerical implementation. Hence, we are essentially performing an eigendecomposition of the covariance matrix $\mathbf{K}$ of $\mathbf{U} = (U(x_1,\omega),\dots,U(x_M,\omega))$ where $\mathbf{K}_{ij} = Cov(U(x_i),U(x_j)) $ for $\{x_i\}_{i=1}^M \subset [0,1]$. In the following, we used an extremely fine mesh with $x_{i+1} - x_i = 0.005$ such that $M = 201$. Further decreasing the mesh size does not alter the conclusions that follow.

In order to approximate the unknown field $A(x,\omega)$, we solve the optimization problem by minimizing the $L^2([0,1])\otimes L^2(\Omega)$ norm, i.e. 
\begin{align} \label{SISC:eq:OptimProblemODE}
\text{min}_{\widetilde{A}} \,\, \int_0^1 E \left [ \left | U(x,\cdot) - \int_0^x \frac{1}{\widetilde{A}(y,\cdot)}\,dy \right | ^2 \right] \, dx = \text{min}_{\widetilde{B}} \,\, \int_0^1 E \left[ \left | U(x,\cdot) - \int_0^x \widetilde{B}(y,\cdot)\, dy\right |^2 \right] \, dx.
\end{align}
We tackle the infinite-dimensional problem to obtain solutions of the form \eqref{SISC:eq:optimalRF} instead of \eqref{SISC:eq:SPGrid} because we do not perform truncation. As a first attempt, we characterize the random field $B$ (instead of $A$) in the spirit of \eqref{SISC:eq:optimalRF} to yield the expression 
\begin{align} \label{SISC:eq:inverseFieldTrialSoln}
\widetilde{B}(x,\omega) = b_0(x) + \sum_{k=1}^M b_k(x) Y_k(\omega)
\end{align}
where $\{b_k(x)\}_{k=0}^{M}$ are to be determined. The quantity $M$ here is not the truncation level according to the decay of the eigenvalues of the covariance function of $U(x,\omega)$ but results from the discrete implementation of the KL expansion as discussed earlier. 


It is observed from the optimization problem on the right-side of \eqref{SISC:eq:OptimProblemODE} that the minimum can be obtained if the unknown deterministic functions $b_k(x)$ satisfy 
\begin{align}\label{SISC:eq:MSderivCond}
\sqrt{\lambda_k} \phi_k(x) = \displaystyle \int_0^x b_k(y) \, dy \,\,\,\, \text{for} \,\,\,\, k=0,\dots,M 
\end{align}
where $\lambda_k, \, \phi_k(x)$ are eigenvalues and eigenfunctions of $Cov(U(x),U(y))$. In other words, the minimizer is achieved if $\widetilde{B}(x,\omega) = \frac{dU(x,\omega)}{dx}$ in the mean square sense \cite{book:Grigoriu2012} which results in $Cov(U(x),U(y)) = \displaystyle \int_0^x \int_0^y Cov(\widetilde{B}(s),\widetilde{B}(t)) \,ds \,dt$ using second moment calculus. This equality is what we would have obtained from the analytic solution for $U$. This implies that by setting $b_k(x)$ as in \eqref{SISC:eq:MSderivCond}, the approach pursued through \eqref{SISC:eq:inverseFieldTrialSoln} recovers the second-order statistics of $B$ and hence that of $A$. 

\begin{figure}[h!]
		\centering
		\includegraphics[width = 0.85\textwidth]		
						{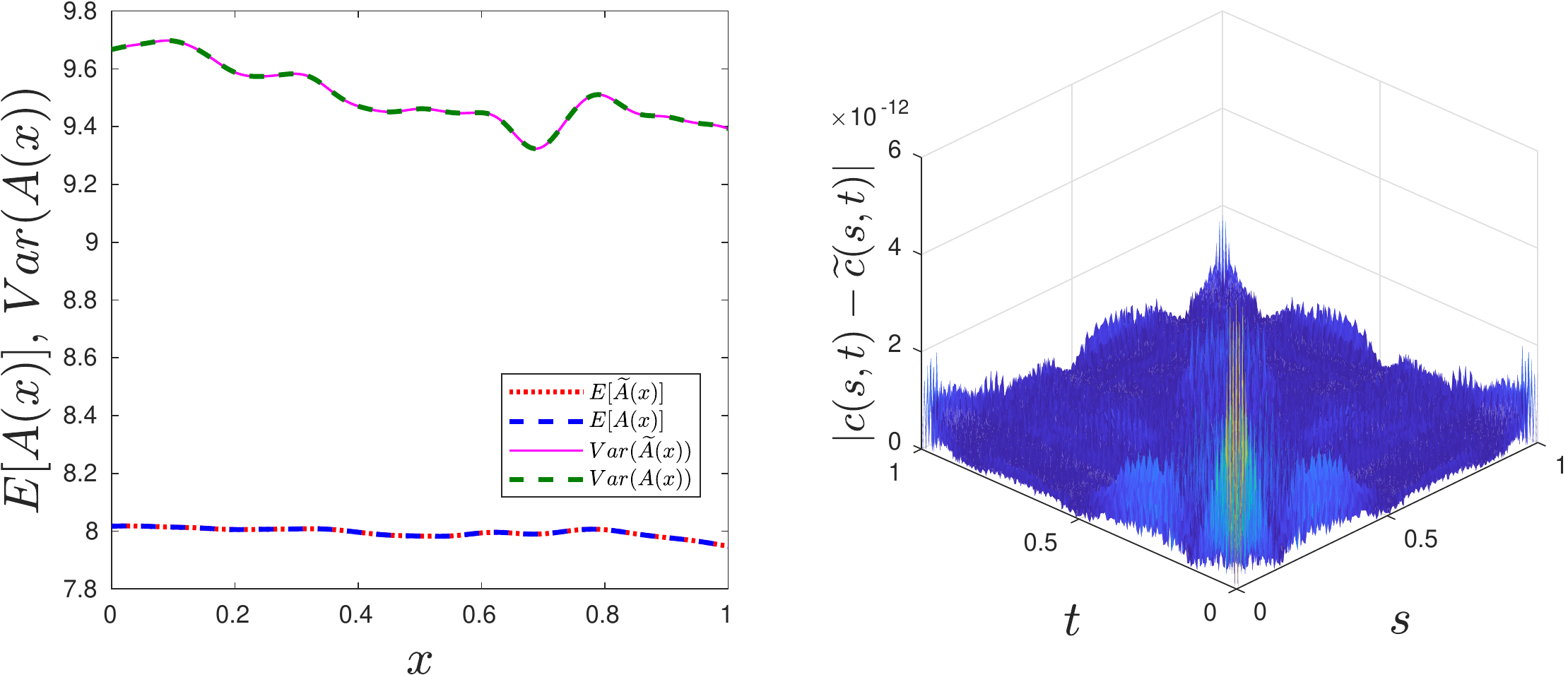}
		\caption{Left panel: Comparison between $E[A(x)]$ (blue dashed line) and $E[\widetilde{A}(x)]$ (red dotted line) and $Var(\widetilde{A}(x))$ (green dashed line) and $Var(\widetilde{A}(x))$ (magenta solid line). Right panel: Plot of the discrepancy between the covariance of $A(x)$ and $\widetilde{A}(x)$, i.e. $|Cov(A(s),A(t)) - Cov(\widetilde{A}(s),\widetilde{A}(t))| = |c(s,t) - \widetilde{c}(s,t)|$. $\widetilde{A}(x)$ is approximated under the first attempt.}
		\label{SISC:fig:InverseFieldStats}
\end{figure}

We verify this numerically by noticing that solving for $\{b_k(x)\}_{k=1}^M$ in \eqref{SISC:eq:OptimProblemODE} results in a matrix least squares problem given the 10000 realizations of $U(x,\omega)$. Figure~\ref{SISC:fig:InverseFieldStats} shows plots of the second-order statistics of the true field $A(x,\omega)$ together with the second-order statistics of the numerical solution to the inverse problem $\widetilde{A}(x,\omega) = (\widetilde{B}(x,\omega))^{-1}$ with condition \eqref{SISC:eq:MSderivCond} on $b_k(x)$. In particular, in the left panel of Figure~\ref{SISC:fig:InverseFieldStats}, the red dotted and blue dashed curves represent $E[\widetilde{A}(x)]$ and  $E[A(x)]$ respectively, while the magenta solid and green dashed curves represent  $Var[\widetilde{A}(x)]$ and  $Var[A(x)]$ respectively. Only 2 curves are visible because the mean and variance of $A(x,\omega)$ and $\widetilde{A}(x,\omega)$  are almost indistinguishable. On the other hand, the right panel shows the absolute value of the difference between $c(s,t) \coloneqq Cov(A(s),A(t))$ and $\widetilde{c}(s,t) \coloneqq Cov(\widetilde{A}(s),\widetilde{A}(t))$ over $s,t \in [0,1]$ which are in good agreement. 

As a second attempt, we pursue a more typical approach in which we characterize the unknown field $A$ directly instead of its inverse, consistent with \cite{paper:BorggaardV2015}. To parameterize $A$, we consider
\begin{align} \label{SISC:eq:coeffFieldTrialSoln}
\widetilde{A}(x,\omega) = a_0(x) + \sum_{k=1}^M a_k(x) Y_k(\omega)
\end{align}
where $\widetilde{A}$ solves the optimization problem \eqref{SISC:eq:OptimProblemODE} and $\{a_k(x)\}_{k=0}^{\infty}$ are to be determined. The unknown field is approximated by solving 
\begin{align}\label{SISC:eq:OptimProblemUpperBound}
\text{min}_{\widetilde{A}} \,\, \int_0^1 E [ |A(x,\cdot) - \widetilde{A}(x,\cdot) |^2 ] \, dx
\end{align}
for $\widetilde{A}(x)$ under the norm $L^2([0,1])\otimes L^2(\Omega)$ where $A(x,\cdot)$ refers to the true field specified in Example~\ref{SISC:ex:vanWyk1}. By using the mean value theorem on $U(x,\cdot) - \displaystyle \int_0^x \frac{1}{\widetilde{A}(y,\cdot)} \, dy$, the boundedness of $A(x,\cdot)$, and assumptions on the boundedness of $\widetilde{A}(x,\cdot)$, minimizing \eqref{SISC:eq:OptimProblemUpperBound} under the constraint that $\widetilde{A}(y,\cdot) > 0$ provides an upper bound for the objective function in \eqref{SISC:eq:OptimProblemODE}. 

\begin{figure}[h!]
		\centering
		\includegraphics[width = 0.85\textwidth]		
						{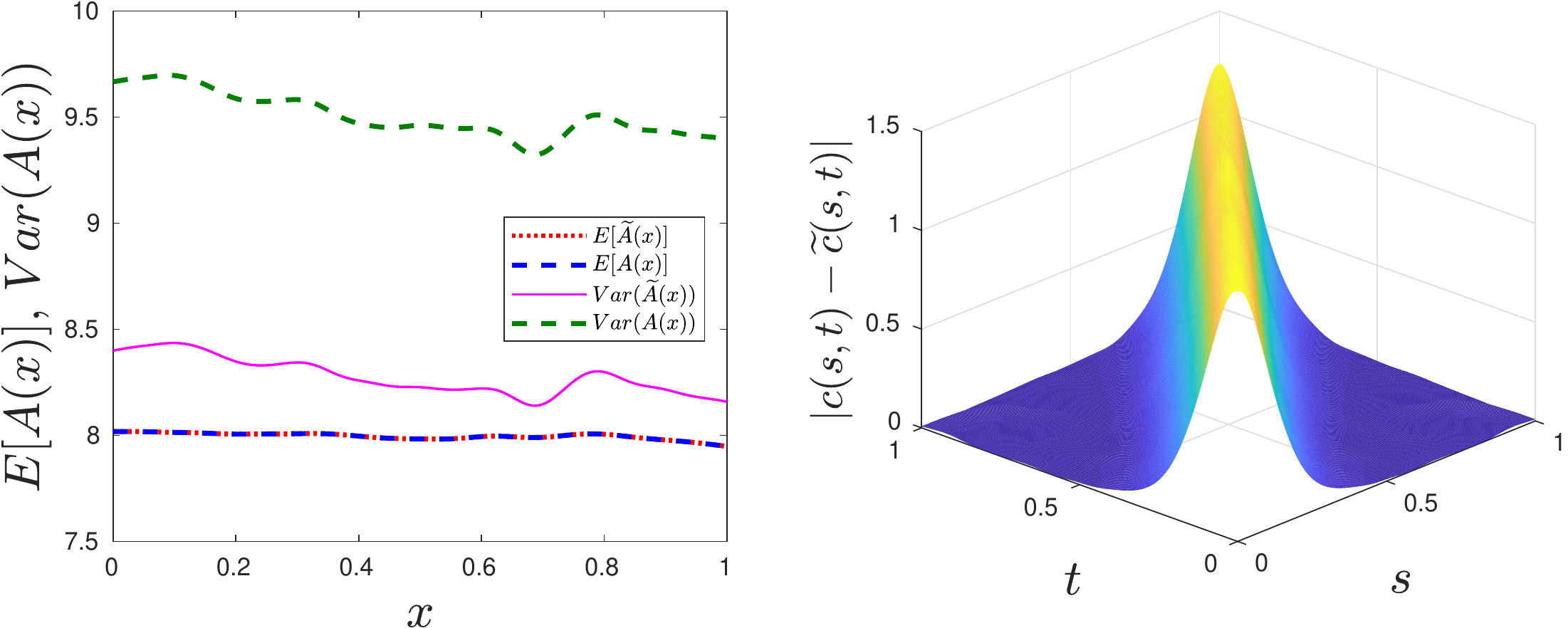}
		\caption{Left panel: Comparison between $E[A(x)]$ (blue dashed line) and $E[\widetilde{A}(x)]$ (red dotted line) and $Var(\widetilde{A}(x))$ (green dashed line) and $Var(\widetilde{A}(x))$ (magenta solid line). Right panel: Plot of the discrepancy between the covariance of $A(x)$ and $\widetilde{A}(x)$, i.e. $|Cov(A(s),A(t)) - Cov(\widetilde{A}(s),\widetilde{A}(t))| = |c(s,t) - \widetilde{c}(s,t)|$. $\widetilde{A}(x)$ is approximated under the second attempt.}
		\label{SISC:fig:CoeffFieldStats}
\end{figure}

Similar to the first attempt, a matrix least squares problem utilizing the 10000 samples of the true $A(x,\omega)$ was solved to optimize \eqref{SISC:eq:OptimProblemUpperBound}.  Figure~\ref{SISC:fig:CoeffFieldStats} shows, with identical legend to Figure~\ref{SISC:fig:InverseFieldStats}, the second order statistics of $\widetilde{A}(x,\omega)$ under the second approach in comparison with that of $A(x,\omega)$. While $E[A(x)]$ and $E[\widetilde{A}(x)]$ are almost similar, there is a considerable difference in the variance and covariance functions between the 2 random fields. Since we minimized \eqref{SISC:eq:OptimProblemUpperBound} instead of the optimization problem on the left-side of \eqref{SISC:eq:OptimProblemODE}, whatever optimal solution we acquire from \eqref{SISC:eq:OptimProblemODE} of the form \eqref{SISC:eq:coeffFieldTrialSoln} cannot do better in matching the statistics of $A(x,\omega)$ than what we have already achieved. The derivations we have presented in the first attempt above do not hold anymore for the current attempt due to the nonlinearity of $U$ with respect to $A$. In addition, increasing the number of samples of $A(x,\omega)$ to 100000 does not alter the results. This underscores that the random variables $\{Y_k\}_{k=1}^{\infty}$ arising from the KL expansion of $U(x,\omega)$ do not necessarily form a stochastic basis for $L^2(\Omega)$. This implies that $\{Y_k\}_{k=1}^{\infty}$ may be inadequate to characterize the unknown $A(x,\omega)$, even when no truncation on the basis of decay of eigenvalues is performed.
However, even when both $U(x,\omega)$ and $A(x,\omega)$ are analytically characterized by the same random variables, the KL expansion of $U$ has to be truncated for numerical implementation. In \cite{paper:BorggaardV2015}, the finite set $\{Y_k\}_{k=1}^M$ resulting from this truncation is then used to characterize $A$. The following example shows that the truncation level for $U$ may not be sufficient for $A$ such that the solution to the inverse problem under this approach underestimates statistics of $A$.

\begin{example} \label{SISC:ex:vanWyk2}
Consider the stochastic ODE for $x \in [0,1], \,\,\omega \in \Omega$:
\begin{align} \label{SISC:eq:ODEexample2}
	\frac{d}{dx} U(x,\omega) = A(x,\omega), \,\,\, U(0,\omega) = 0
\end{align}
whose analytical solution is given by $U(x,\omega) = \displaystyle \int_0^x A(y,\omega) \,dy$. Set $A(x,\omega) = G(x,\omega)$ where $G(x,\omega)$ is a zero mean, unit variance stationary Gaussian process with spectral density \cite[p. 196]{book:SoongG1993} $s_G(\nu) = \frac{F}{(\nu^2 - \nu_0^2)^2 + (2\zeta \nu \nu_0^2)^2}$ in which $\nu_0 = 20$, $\zeta = 0.1$, and $F$ is a scaling factor such that the correlation function of $A$, $r_{G}(\tau)$, satisfies $r_G(\tau) = \displaystyle \int_{-\infty}^{\infty} e^{i \nu \tau}s_G(\nu)\,d \nu = 1.$ To formulate the inverse problem, we approximate $A(x,\omega)$ using observed samples of $U(x,\omega)$.  
\end{example} 

As $A(x,\omega)$ is a Gaussian process, the response $U(x,\omega)$ is also a Gaussian process with mean $m_U = \displaystyle \int_0^x E[A(y)] \,dy = 0$ and correlation function $r_U(x,y) = \displaystyle \int_0^x \int_0^y r_G(s-t) \,ds\,dt$ \cite{book:Grigoriu2002}. Consequently, the KL expansion of both $U(x,\omega)$ and $A(x,\omega)$ can be expressed as
\begin{align*}
A(x,\omega) = \sum_{k=1}^{\infty} \sqrt{\lambda^G_k} \phi^G_k(x) Y_k(\omega), \,\,\,\, U(x,\omega) = \sum_{k=1}^{\infty} \sqrt{\lambda^U_k} \phi^U_k(x) \widetilde{Y}_k(\omega)
\end{align*}
where $Y_k(\omega), \widetilde{Y}_k(\omega) \sim N(0,1) \,\, \forall k$ and $\lambda^G_k, \, \phi^G_k(x)$ and $\lambda^U_k, \, \phi^U_k(x)$ are eigenvalues and eigenfunctions of $r_G$ and $r_U$, respectively. Thus, $A$ and $U$ are characterized by random variables with the same law unlike in Example~\ref{SISC:ex:vanWyk1}. In the solution approach in \cite{paper:BorggaardV2015}, the truncation level of the KL expansion of $U$ is used to determine the number of random variables to characterize $A$. The truncation level $M$ of the KL expansion is usually deduced from the total variance of the random field which is given by $\displaystyle \int_0^1 Var(U(x)) \, dx = \sum_{k=1}^{\infty} \lambda_k^U$. The value of $M$ is then chosen to be the smallest integer such that $\frac{\sum_{k=1}^M \lambda_k^U}{\sum_{k=1}^\infty \lambda_k^U} \ge \alpha$ with $\alpha$ being close to 1.

\begin{figure}[h!]
		\centering
		\includegraphics[width = 0.8\textwidth]		
						{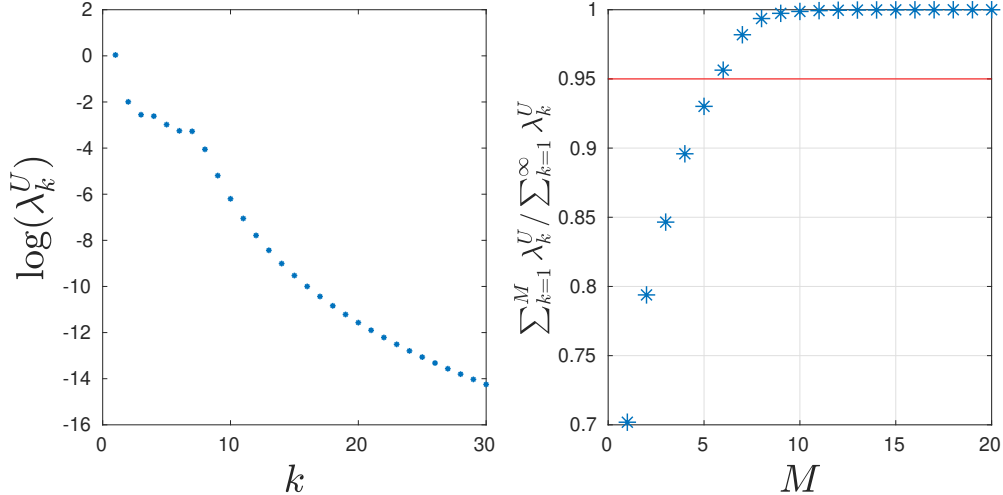}
		\caption{Left: First 30 eigenvalues of $r_U$. Right: Truncation level criterion $\frac{\sum_{k=1}^M \lambda_k^U}{\sum_{k=1}^\infty \lambda_k^U}$ vs $M$.}
		\label{SISC:fig:respCovarEigVals}
\end{figure}

The above principle is applied to choose $M$ using $\alpha = 0.95$. To simulate solving the inverse problem, we do not generate samples of $U$ through samples of $A$ in \eqref{SISC:eq:ODEexample2} and estimate $r_U(x,y)$ from samples of $U$; rather, $r_U$ is obtained from the relationship between $r_U$ and $r_G$ above. Figure~\ref{SISC:fig:respCovarEigVals} exhibits the behavior of the eigenvalues of $r_U$. The left panel shows the first 30 eigenvalues $\lambda_k^U$ while the right panel displays $\frac{\sum_{k=1}^M \lambda_k^U}{\sum_{k=1}^\infty \lambda_k^U}$ as a function of $M$. It is evident that the truncation level for the KL expansion of $U$ according to the above procedure is $M=6$.

\begin{figure}[h!]
		\centering
		\includegraphics[width = 0.8\textwidth]		
						{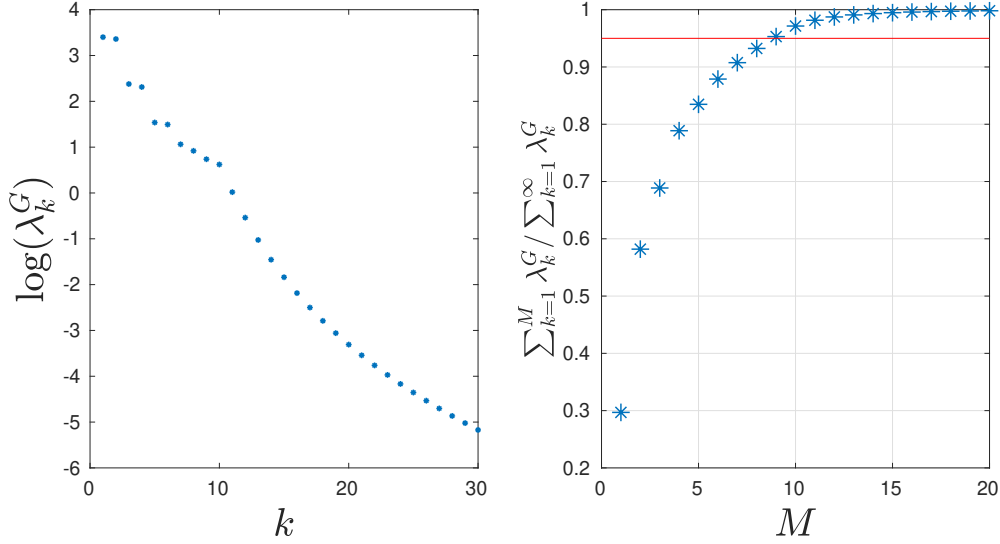}
		\caption{Left: First 30 eigenvalues of $r_G$. Right: Truncation level criterion $\frac{\sum_{k=1}^M \lambda_k^G}{\sum_{k=1}^\infty \lambda_k^G}$ vs $M$.}
		\label{SISC:fig:coeffFieldCovarEigVals}
\end{figure}

Figure~\ref{SISC:fig:coeffFieldCovarEigVals} displays the behavior of the eigenvalues $\lambda_k^G$ of $r_G$ with the same legend as in Figure~\ref{SISC:fig:respCovarEigVals}. In this case, the truncation level for the KL expansion of $A$ is $M=9$. Note that the eigenvalues of $r_G$ decay slower than that of $r_U$; intuitively, this is because large values of the frequency $\nu$ are required to capture the total energy of the spectral density $\displaystyle 2\int_0^{\infty} s_G(\nu) \, d\nu$. As $U(x,\omega)$ is obtained by integrating $A(x,\omega)$, the variation in $A(x,\omega)$ is diminished which yields a faster decay of eigenvalues for $r_U$.  

\begin{figure}[h!]
		\centering
		\includegraphics[width = 0.9\textwidth]		
						{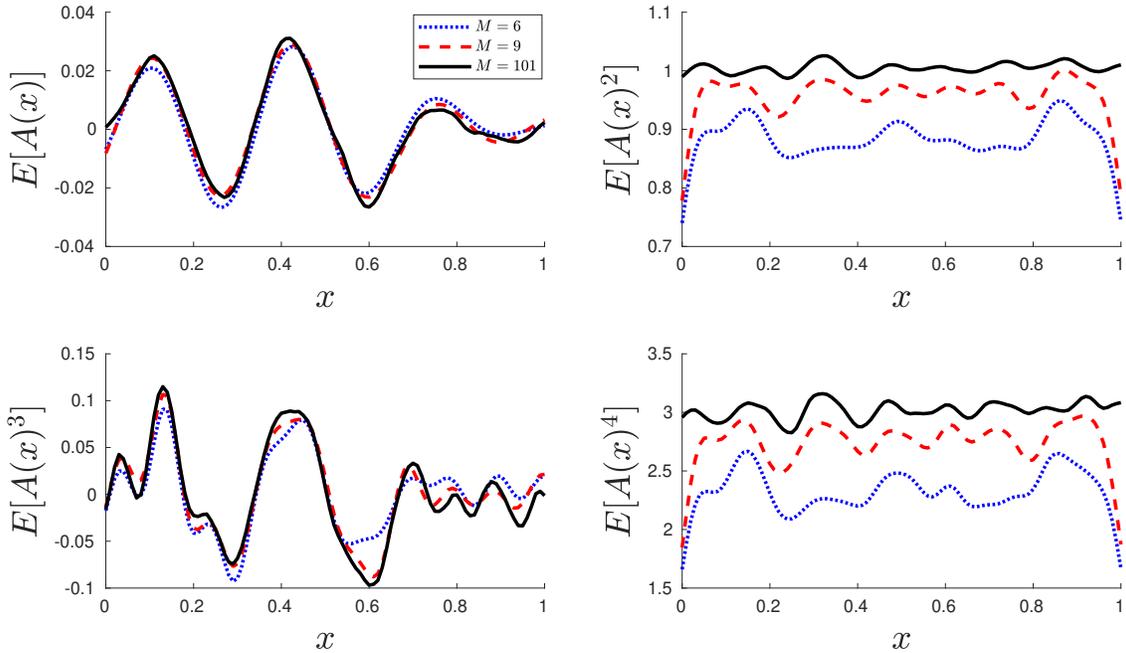}
		\caption{$p$-th order moments of $\sum_{k=1}^M \sqrt{\lambda_k^G} \phi_k^G(x) Y_k(\omega)$ for $p=1,\dots,4$ and $M=6,$ (blue dotted line), $9,$ (red dashed line), and $101$ (black solid line).}
		\label{SISC:fig:TruncatedFieldStats}
\end{figure}

Hence, if the truncation level of the KL expansion of $U$ is used to characterize $A$, the inverse problem solution would be $\widetilde{A}(x,\omega) = \sum_{k=1}^6 \sqrt{\lambda_k^G} \phi_k^G(x) Y_k(\omega)$. Based on the discussion above, $\widetilde{A}$ may not be sufficient to capture the statistics of $A$ and we confirm this in Figure~\ref{SISC:fig:TruncatedFieldStats}. The four subplots in this figure represent the first 4 moments of $\sum_{k=1}^M \sqrt{\lambda_k^G} \phi_k^G(x) Y_k(\omega)$ for $M=6,9,101$ with $M = 101$ being a sufficient approximation for $M = \infty$. We notice that for $M=6$, $\widetilde{A}(x,\omega)$ underestimates the statistics of $A(x,\omega)$, especially for moments of even order. The underestimation can be avoided provided appropriate selection of the truncation level. Perhaps deducing the truncation level for $A$ using other information than the observations $U$ can resolve this issue.

To summarize, Examples~\ref{SISC:ex:vanWyk1} and~\ref{SISC:ex:vanWyk2} demonstrate that in parameterizing an unknown random field $A(x,\omega)$ by random variables $\{Y_k\}_{k=1}^M$ obtained from the response $U(x,\omega)$, the value chosen for $M$ and the probability law chosen for $Y_k$ can significantly affect the accuracy of the approximation for $A(x,\omega)$. A strategy to tackle the ill-posedness could instead require additional information being specified on the dimension of the random variables parameterizing $A$ and that they belong to a family of distributions subject to unknown parameters. This will be elaborated next.

\section{Required additional information on the unknown random quantity} \label{SISC:sec:AddlInfo}

Section~\ref{SISC:sec:AbsentInfo} showed that the existing methods in general cannot recover the true probability law of a random vector $Z$ if the information is limited to its bounded domain. We therefore devote this section to addressing the objective of this work: identify realistic additional information on $Z$ that is required to characterize its law. The mapping in Example~\ref{SISC:ex:butler_estep}, $Q(Z)=Z_1\cdot Z_2$, is revisited in which the true law is set as $Z_1,Z_2 \sim U(0,1)$, independent. All pdfs from this section onwards are constructed with respect to the Lebesgue measure and are denoted by $f$.


It is impossible to devise a general method to determine the minimum amount of additional information that is required on $Z$. Different applications possess forward models with unique properties and specific information on the unknown. Hence, we outline a few scenarios with their corresponding solution methodologies. These are categorized based on what is known about $Z$: moment information (Section~\ref{SISC:subsec:GivenMomentofZ}) or the family of distribution to which it belongs (Section~\ref{SISC:subsec:GivenFamilyofZ}). In each category, further subcategories are considered depending on the given information on $Q(Z)$:
pdf (Sections~\ref{SISC:subsubsec:entropyPushFwd} and~\ref{SISC:subsubsec:KnownPdfPushFwd}) or samples (Sections~\ref{SISC:subsubsec:entropyBayes} and~\ref{SISC:subsubsec:KnownPdfBayes}).


\subsection{Information on moments of $Z$} \label{SISC:subsec:GivenMomentofZ}

If information about moments of the random vector $Z$ is available, the principle of maximum entropy \cite{book:CoverT2006} can be employed to determine the pdf of $Z$, assuming that this philosophy is accepted and that it is believed that the solution to the inverse problem resides in the subspace of maximum entropy pdfs. The principle of maximum entropy constructs the pdf $f_Z(z)$ of $Z$ by solving

\begin{equation} \label{SISC:eq:entropyProblem}
	\begin{aligned}
		& \underset{f_Z}{\text{minimize}}
		& & \int_{\Gamma} f_Z(z) \, \log(f_Z(z)) \,dz \\
		& \text{subject to}
		& & \int_{\Gamma} g_k(z) f_Z(z) \,dz = \mu_k, \,\,k=1,\dots,N, \\
		&&& \int_{\Gamma} f_Z(z) \,dz = 1.
	\end{aligned}
\end{equation}
for some functions $g_k$, whose solution is derived as
\begin{align} \label{SISC:eq:entropyPdfOptimal}
f_Z(z) = \frac{1}{\int_{\Gamma} \exp[\lambda_1 g_1(z) + \dots + \lambda_N g_N(z)] \, dz}\exp[\lambda_1 g_1(z) + \dots + \lambda_N g_N(z)],\,\,\, z \in \Gamma
\end{align}
where the Lagrange multipliers $\lambda_k$ satisfy the relationship $\displaystyle \frac{\partial}{\partial \lambda_k} \int_{\Gamma} \exp[\lambda_1 g_1(z) + \dots + \lambda_N g_N(z)] \, dz = \mu_k$ for $k = 1,\dots,N$. If these multipliers exist, it can be shown that \eqref{SISC:eq:entropyPdfOptimal} is the unique minimizer satisfying the above constraints \cite{book:CoverT2006}. The succeeding sections detail how \eqref{SISC:eq:entropyPdfOptimal} can be used to solve the inverse problem given the pdf of $Q(Z)$ (Section~\ref{SISC:subsubsec:entropyPushFwd}) or samples of $Q(Z)$ (Section~\ref{SISC:subsubsec:entropyBayes}).

\subsubsection{Pdf of $Q(Z)$} \label{SISC:subsubsec:entropyPushFwd}

Assume that the pdf $f_Q(q) = -\log(q), \, q \in (0,1]$ of $Q$ were known. The example below illustrates the above construction.

\begin{example}
Suppose that the only information known about $Z$ aside from $Z \in \Gamma$ are that 1) $Z_1$ and $Z_2$ are independent and that 2) the first-order moments of $Z_1$ and $Z_2$ are within a certain range, i.e. $E[Z_1] = \mu_1 \in [0,0.75], E[Z_2] = \mu_2 \in [0.4,1].$ It is demonstrated how the inverse problem can be solved using an entropy-based pdf.
\end{example} 

For a fixed value of $(\mu_1,\mu_2)$, the principle of maximum entropy yields the conditional pdf of $Z$ as

\begin{align}\label{SISC:eq:entropyPdf}
f_Z(z_1,z_2 | \mu_1, \mu_2) = \frac{\lambda_1}{(e^{\lambda_1}-1)} \frac{\lambda_2}{(e^{\lambda_2}-1)} \exp[\lambda_1 z_1 + \lambda_2 z_2]
\end{align} 
where there is a bijective relationship between $\mu_i$ and $\lambda_i$ for $i=1,2$ via $\mu_i = \frac{1}{1-e^{-\lambda_i}} - \frac{1}{\lambda_i}$, $\mu_i \in (0,1)$. To estimate $(\mu_1,\mu_2)$, an optimization problem can be solved which measures the discrepancy between the given pdf $f_Q$ of $Q$ and the one obtained by propagating the pdf in \eqref{SISC:eq:entropyPdf} through the forward model $Q$ which we denote by $\widetilde{f}_Q(\cdot| \mu_1,\mu_2)$. Mathematically, this is expressed as 
\begin{align}\label{SISC:eq:optimProblemPushFwd}
(\mu_1,\mu_2) = \underset{(\mu_1,\mu_2)}{\text{argmin}} \,\,\, d(f_Q(q),\widetilde{f}_Q(q| \mu_1,\mu_2))
\end{align}
for some distance function $d$ such as the $L^p$ error, Kullback-Leibler divergence, etc. Let $F_{Z_i}(\cdot | \mu_1,\mu_2)$ and $f_{Z_i}(\cdot | \mu_1,\mu_2)$ represent the marginal cdf and pdf, respectively, of $Z_i$ for $i=1,2$ based on \eqref{SISC:eq:entropyPdf}. Elementary calculations similar to \eqref{SISC:eq:prodPdfDerivation} yield
\begin{align} \label{SISC:eq:ComputePDFofProdGeneral}
\widetilde{f}_Q(q| \mu_1,\mu_2) = \displaystyle \frac{d}{dq} \int_0^1 F_{Z_2}\left(\frac{q}{z_1}\big\rvert \mu_1,\mu_2 \right)f_{Z_1}(z_1 | \mu_1,\mu_2) \,dz_1 = \int_q^1 f_{Z_2}\left(\frac{q}{z_1}\big\rvert \mu_1,\mu_2 \right) f_{Z_1}(z_1 | \mu_1,\mu_2) \,dz_1 
\end{align}
for $q \in [0,1]$. Figure~\ref{SISC:fig:EntropyPushFwdError} displays the logarithm of the $L^1$ error $\|f_Q(q) - \widetilde{f}_Q(q|\mu_1,\mu_2)\|_{L^1}$ for $(\mu_1,\mu_2)$ in the specified ranges above. This discrepancy is minimized at $(\mu_1,\mu_2) = (0.5,0.5)$ which corresponds to $(\lambda_1,\lambda_2) = (0,0)$, thereby recovering the true pdf of $Z$: $f_Z(z_1,z_2 | \mu_1,\mu_2) = \mathbbm{1}_{(z_1,z_2) \in \Gamma}$ with $\mathbbm{1}$ being the indicator function.

For this simple example, the pdf of $Q$ given the pdf of $Z$ can be computed analytically. In general, however, if $f_Z(\cdot|\theta)$ for $\theta \in \Theta$ represents the pdf of $Z$, propagating this through more complicated forward models $Q$ implies evaluation of $Q$ multiple times. An approach to ameliorate this computational burden is to use a surrogate model \cite{book:Grigoriu2012} for $Q$ as a function of $Z$. This can be supplemented by the following procedure if $\Theta$ is low-dimensional, as done in \cite{paper:EmeryGF2016}:

\begin{itemize}
	\item Select $M$ points $\{\theta_i\}_{i=1}^M \subset \Theta$.
	\item For each $\theta_i$, $i=1,\dots,M$, propagate $f_Z(\cdot | \theta_i)$ through   $Q$ to approximate the pdf $\widetilde{f}_Q(\cdot | \theta_i)$ of $Q$.
	\item Using $\{\widetilde{f}_Q(\cdot| \theta_i)\}_{i=1}^M$, construct an interpolant for $\widetilde{f}_{Q}(\cdot|\theta)$ over $\Theta$.
\end{itemize}

\subsubsection{Samples of $Q(Z)$} \label{SISC:subsubsec:entropyBayes}

In contrast to the previous section, consider that the available information on the quantity of interest $Q$ is its $N_s$ samples represented by $\{q^i\}_{i=1}^{N_s}$ instead of the pdf of $Q$. Such information is what is typically encountered in practical applications of stochastic inverse problems \cite{paper:EmeryGF2016,paper:NarayananZ2004}.
This naturally leads to employing the Bayesian framework \cite{book:KaipioS2005} to solve the inverse problem. If information is available on moments of $Z$ in the form of a prior pdf, the principle of maximum entropy can be utilized to construct the likelihood function in Bayes' theorem as the next example elaborates.

\begin{example} \label{SISC:ex:EntropyBayes}
We postulate that only the following information is known about $Z$: 1) $Z \in \Gamma$, 2) $Z_1,Z_2$ are independent, and that 3) $(\mu_1,\mu_2) \coloneqq (E[Z_1],E[Z_2])$ is equally likely to take any value in the range $[0.25,0.75]^2$. Bayes' theorem coupled with the principle of maximum entropy provide an approach to address the inverse problem.
\end{example}

Knowledge that $(\mu_1,\mu_2)$ is equally likely in $[0.25,0.75]^2$ translates to a prior pdf on $(\mu_1,\mu_2)$ denoted by $f_{\mu}^{prior}(\mu_1,\mu_2) = \mathbbm{1}_{(\mu_1,\mu_2)\in [0.25,0.75]^2} \frac{1}{0.5^2}$. By Bayes' theorem, the posterior pdf on $(\mu_1,\mu_2)$ is
\begin{align}\label{SISC:eq:EntropyPosterior}
f_{\mu}^{post}(\mu_1,\mu_2 | \{q^i\}_{i=1}^{N_s}) = \displaystyle \frac{\ell(\mu_1,\mu_2 |  \{q^i\}_{i=1}^{N_s})f_{\mu}^{prior}(\mu_1,\mu_2)}{\int_{[0.25,0.75]^2} \ell(\mu_1,\mu_2 |  \{q^i\}_{i=1}^{N_s})f_{\mu}^{prior}(\mu_1,\mu_2) \,d\mu_1 d\mu_2}
\end{align}
in which $\ell(\mu_1,\mu_2 |  \{q^i\}_{i=1}^{N_s})$ symbolizes the likelihood function. Given $(\mu_1,\mu_2)$, let $f_Z(\cdot|\mu_1,\mu_2)$ as in \eqref{SISC:eq:entropyPdf} be the entropy-based pdf of $Z$ while $\widetilde{f}_{Q}(\cdot|\mu_1,\mu_2)$ as in \eqref{SISC:eq:ComputePDFofProdGeneral} be the resulting pdf when the entropy-based pdf is propagated through $Q$. The likelihood function is then established as
\begin{align}\label{SISC:eq:EntropyLikelihood}
\ell(\mu_1,\mu_2 | \{q^i\}_{i=1}^{N_s}) = \prod_{i=1}^{N_s} \widetilde{f}_Q(q^i| \mu_1,\mu_2)
\end{align}
by the independence of the samples $\{q^i\}_{i=1}^{N_s}$.

Figure~\ref{SISC:fig:EntropyPosterior} displays the posterior pdf $f_{\mu}^{post}(\cdot| \{q^i\}_{i=1}^{N_s})$ on $(\mu_1,\mu_2)$ using $N_s = 100$ samples of $Q$. The maximum a posteriori (MAP) estimate is around $(\mu_1,\mu_2) = (0.5,0.5)$ which recovers the true pdf of $Z$. We remark that the accuracy of this approach is in part influenced by the observed number of samples $N_s$ of $Q$.


Does the entropy-based pdf for $Z$ guarantee that the true pdf of $Z$ can be recovered, or that the optimization problem formulated or the posterior density have a unique global minimum? The answer depends on the specific application. If not, additional information of the moments of the unknown $Z$ need to be specified to obtain a solution to the inverse problem that can be used for prediction as elaborated in Section~\ref{SISC:subsec:validation}. Sufficient conditions exist which impose criteria that the moments of $Z$ have to satisfy to uniquely determine its distribution; see \cite[Theorem 3.3.11]{book:Durrett2010} for an example.

\subsection{Parametric family of distributions of $Z$} \label{SISC:subsec:GivenFamilyofZ}

Another strategy to combat ill-posedness is to require that the practitioner has information about the family of distributions in which the law of $Z$ resides, subject to unknown parameters $\theta$. This information is represented as $f_Z(\cdot | \theta)$ in which the functional form of the pdf of $Z$ is known. We illustrate how this can be employed to solve the inverse problem if the pdf of $Q(Z)$ (Section~\ref{SISC:subsubsec:KnownPdfPushFwd}) or samples of $Q(Z)$ (Section~\ref{SISC:subsubsec:KnownPdfBayes}) is given.

\subsubsection{Pdf of $Q(Z)$} \label{SISC:subsubsec:KnownPdfPushFwd}

The next example expounds on the above idea assuming that the pdf $f_Q(q) = -\log(q), q \in (0,1]$ of $Q$ were known.

\begin{example}
Suppose that it is known that $Z_1,Z_2$ are independent with $Z_1\sim Beta(\nu_1,\nu_1)$ and $Z_2 \sim Beta(\nu_2,\nu_2)$ whose joint pdf is denoted by $f_Z(\cdot | \nu_1,\nu_2)$. The objective is to estimate $(\nu_1,\nu_2)$ such that the pdf $f_Q$ of $Q$ matches the pdf obtained by propagating $f_Z(\cdot | \nu_1,\nu_2)$ through $Q$.
\end{example}

\begin{figure}[!h]
\centering
\begin{subfigure}{.5\textwidth}
  \centering
  \includegraphics[width=.9\linewidth]{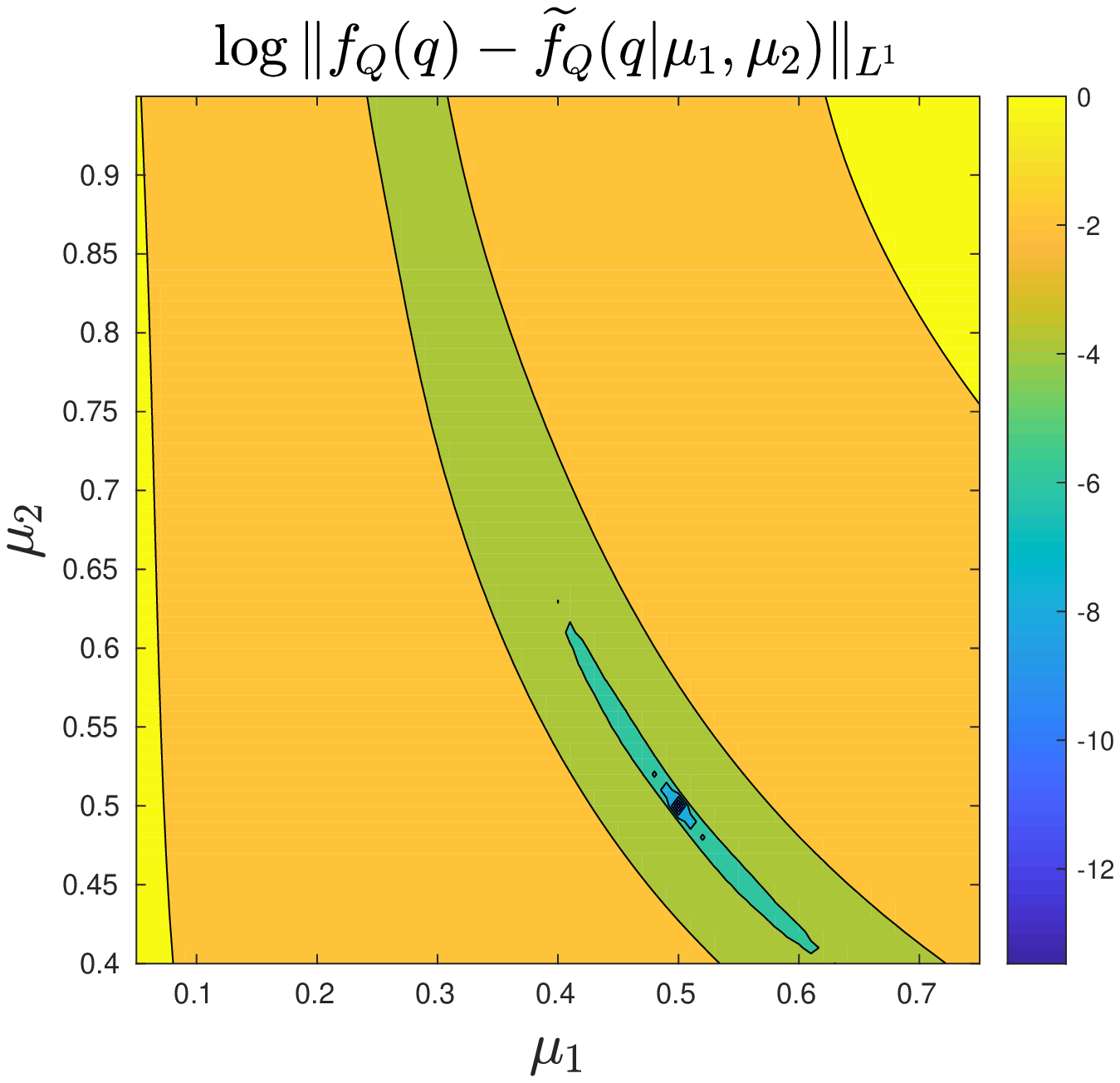}
  \caption{Entropy-based pdf}
  \label{SISC:fig:EntropyPushFwdError}
\end{subfigure}%
\begin{subfigure}{.5\textwidth}
  \centering
  \includegraphics[width=0.9\linewidth]{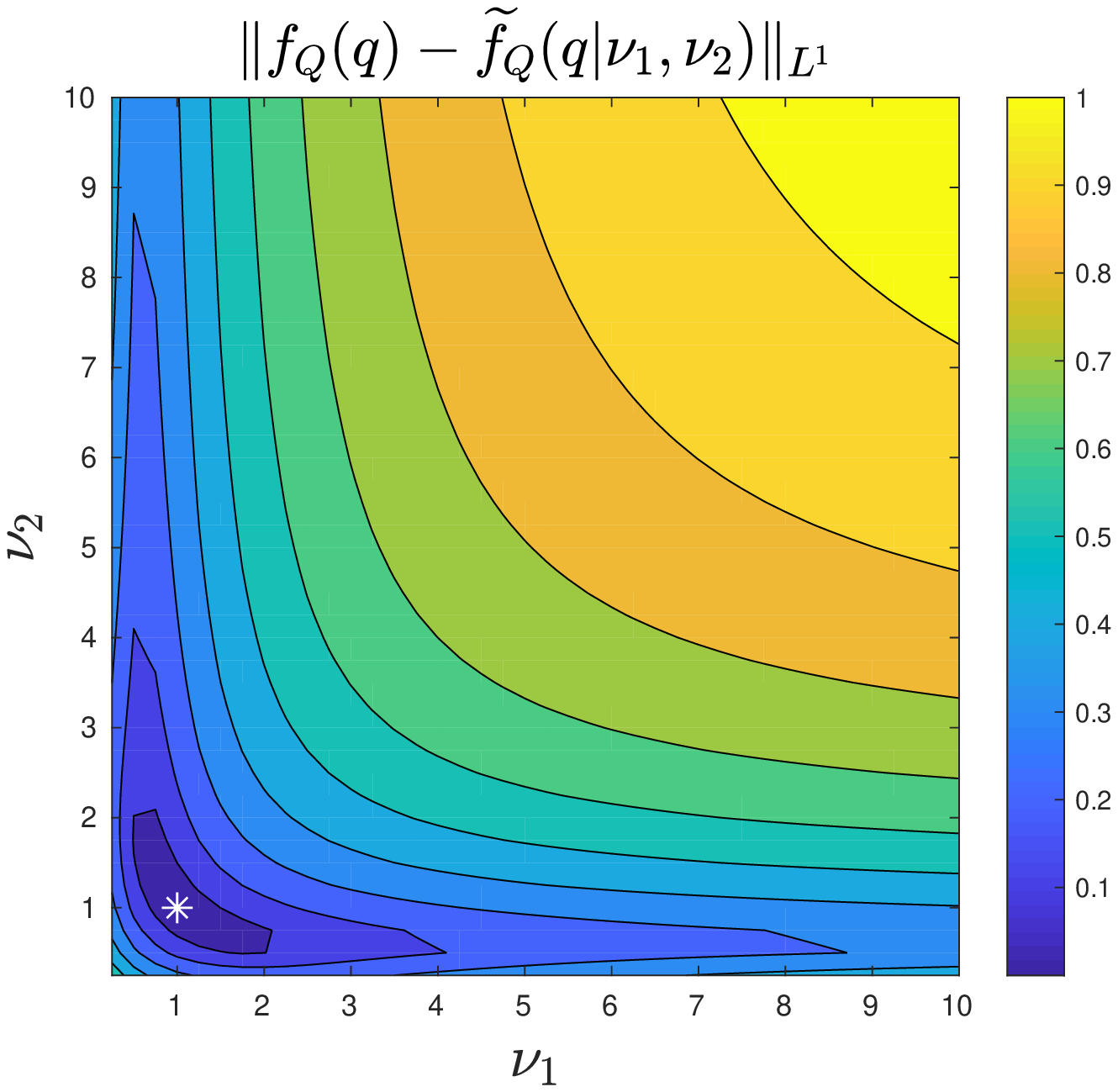}
   \caption{Known distribution, unknown parameters}
   \label{SISC:fig:KnownModelPushFwdError}
\end{subfigure}
\caption{Discrepancy between the given pdf $f_Q$ of $Q$ and the pdf obtained by propagating the pdf $f_Z(\cdot|\theta)$ of $Z$ through $Q$. Left panel: $f_Z(\cdot|\theta)$ is obtained through the principle of maximum entropy. Right panel: $f_Z(\cdot|\theta)$ is a specified distribution subject to unknown parameters. The white asterisk denotes the location of the global minimum.}
\end{figure}

The solution methodology for this approach is identical to the optimization procedure in \eqref{SISC:eq:optimProblemPushFwd} in which we write the propagated pdf of $f_Z(\cdot | \nu_1,\nu_2)$ through $Q$ as $\widetilde{f}_Q(\cdot|\nu_1,\nu_2)$. By nature of the beta distribution, $(\nu_1,\nu_2) \in (0,\infty)^2$. Figure~\ref{SISC:fig:KnownModelPushFwdError} displays the $L^1$ error $\|f_Q(q) - \widetilde{f}_Q(q|\nu_1,\nu_2)\|_{L^1}$ for $(\nu_1,\nu_2) \in (0,10]^2$ where the global minimum in this domain is attained at $(\nu_1,\nu_2) = (1,1)$, thereby recovering the true pdf of $Z$. Heuristic arguments can be made to deduce that 
no other values for $(\nu_1,\nu_2)$ outside  $(0,10]^2$ yield the global minimum. If $\nu_1,\nu_2$ are simultaneously large, $Z_1,Z_2 \rightarrow \frac{1}{2}$ a.s. which implies that $Q \rightarrow \frac{1}{4}$ a.s. Likewise, if only $\nu_1$ is large then $Q \rightarrow \frac{1}{2}Z_2$ a.s. whose pdf does not match $f_Q(q) = -\log(q)$, and vice versa. 

\subsubsection{Samples of $Q(Z)$} \label{SISC:subsubsec:KnownPdfBayes}

In contrast to the previous example, suppose instead that the available information on $Q$ pertains to its $N_s$ samples $\{q^i\}_{i=1}^{N_s}$. Information on the family of distributions enables the construction of the likelihood function that is required to find the posterior pdf on the unknown parameters using Bayes' theorem. The next example highlights this idea of standard Bayesian inversion.

\begin{example}\label{SISC:ex:knownPdfBayes}
Consider a model in which the following information is at the practitioner's disposal: 1) $Z_1,Z_2$ are independent, 2) $Z_1 \sim Beta(1,\nu_1), Z_2 \sim Beta(1,\nu_2)$, and 3) the prior pdf on $(\nu_1,\nu_2)$ is characterized by $f_{\nu}^{prior}(\nu_1,\nu_2) = \mathbbm{1}_{(\nu_1,\nu_2)\in [\frac{1}{3},3]^2} \frac{1}{0.5^2} \frac{1}{(\nu_1 + 1)^2} \frac{1}{(\nu_2 + 1)^2}$. The posterior pdf on $(\nu_1,\nu_2)$ results directly from Bayes' theorem.
\end{example} 

\begin{figure}
        \centering
        \begin{subfigure}[b]{0.45\textwidth}
            \centering
            \includegraphics[width=\textwidth]{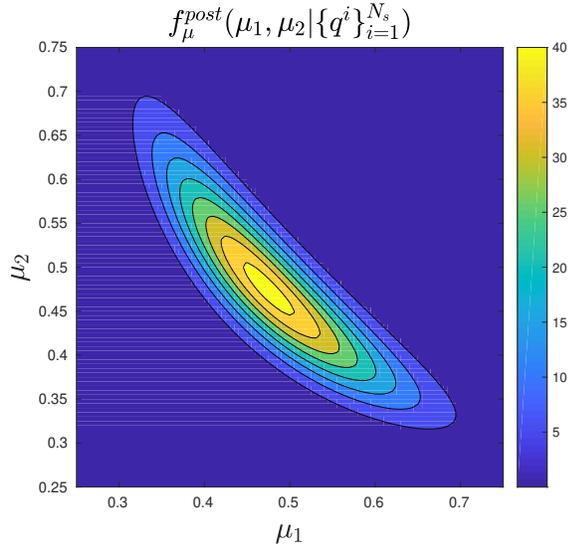}
            \caption{Entropy-based likelihood}   
            \label{SISC:fig:EntropyPosterior}
        \end{subfigure}
       	\hfill
        \begin{subfigure}[b]{0.45\textwidth}  
            \centering 
            \includegraphics[width=\textwidth]{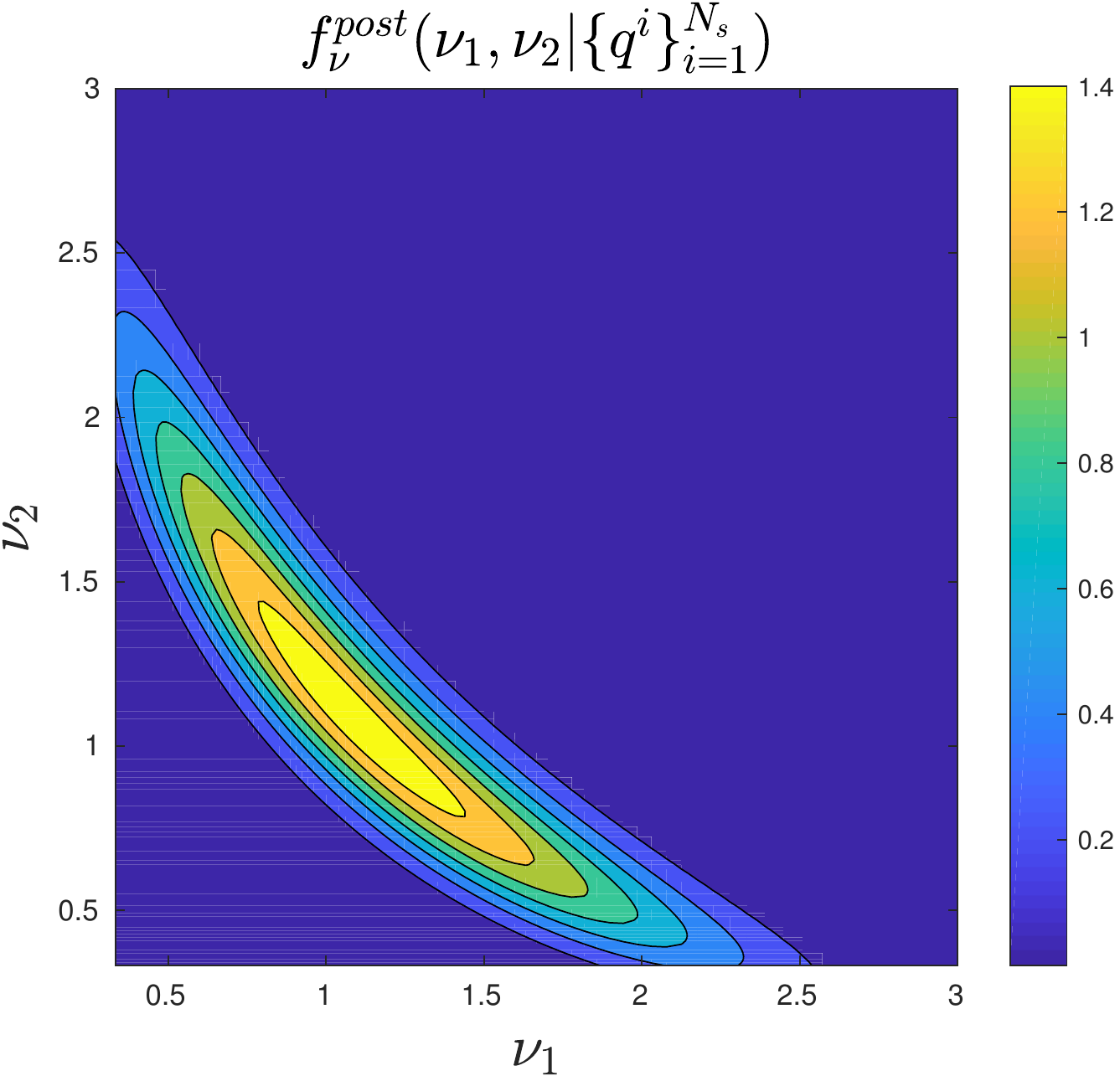}
            \caption{Known distribution, unknown parameters}    
            \label{SISC:fig:KnownPdfPosterior}
        \end{subfigure}
        \vskip\baselineskip
        \begin{subfigure}[b]{0.45\textwidth}   
            \centering 
            \includegraphics[width=\textwidth]{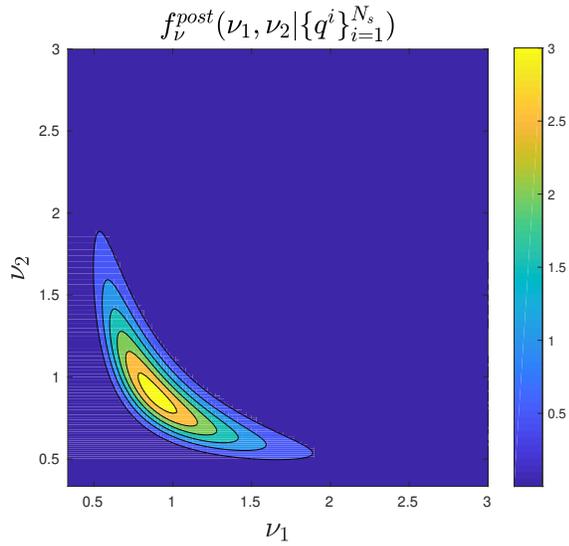}
            \caption{Entropy-based likelihood}    
            \label{SISC:fig:EntropyPosteriorChangeVar}
        \end{subfigure}
        \hfill
        \begin{subfigure}[b]{0.45\textwidth}   
            \centering 
            \includegraphics[width=\textwidth]{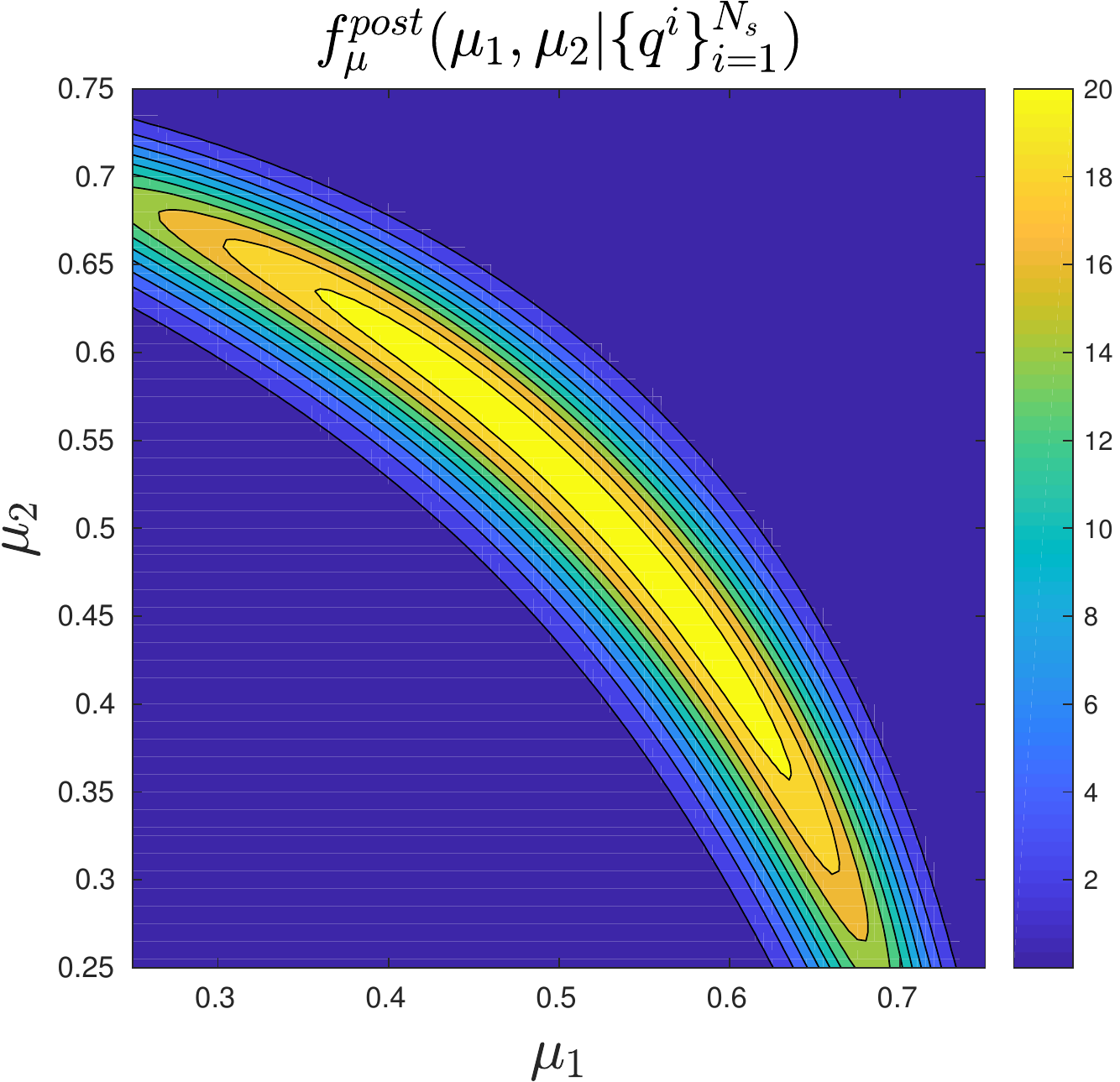}
            \caption{Known distribution, unknown parameters}    
            \label{SISC:fig:KnownPdfPosteriorChangeVar}
        \end{subfigure}
        \caption{Posterior distribution of the parameters in Examples~\ref{SISC:ex:EntropyBayes} and~\ref{SISC:ex:knownPdfBayes}. Plots (a) and (c): posterior density in the $(\mu_1,\mu_2)$ parameter space and in the $(\nu_1,\nu_2) = (\frac{1}{\mu_1}-1,\frac{1}{\mu_2}-1)$ parameter space, respectively, in which the likelihood is constructed using the principle of maximum entropy as in Section~\ref{SISC:subsubsec:entropyBayes}. Plots (b) and (d):  posterior density in the $(\nu_1,\nu_2)$ parameter space and in the $(\mu_1,\mu_2) = (\frac{1}{1+ \nu_1},\frac{1}{1+\nu_2})$ parameter space, respectively, in which the likelihood is constructed using the known family of distributions as in Section~\ref{SISC:subsubsec:KnownPdfBayes}. } 
    \end{figure}

Since $E[Z_i] = \frac{1}{1 + \nu_i}$ for $i=1,2$, the prior pdf $f_{\nu}^{prior}$ translates to a uniform prior pdf on $(E[Z_1],E[Z_2])$ with values in the range $[0.25,0.75]^2$. The construction of the likelihood function and the posterior pdf $f_{\nu}^{post}(\nu_1,\nu_2 | \{q^i\}_{i=1}^{N_s})$ is identical to that in \eqref{SISC:eq:EntropyLikelihood} and \eqref{SISC:eq:EntropyPosterior}, respectively, wherein the specified pdf on $Z$ conditioned on $(\nu_1,\nu_2)$ is propagated through $Q$ to obtain $\widetilde{f}_Q(\cdot| \nu_1,\nu_2)$. Figure~\ref{SISC:fig:KnownPdfPosterior} exhibits $f_{\nu}^{post}(\nu_1,\nu_2 | \{q^i\}_{i=1}^{N_s})$ using the same $N_s = 100$ samples of $Q$ generated in Section~\ref{SISC:subsubsec:entropyBayes}. The MAP estimate hovers close to $(\nu_1,\nu_2) = (1,1)$; with more samples of $Q$, the contours of $f_{\nu}^{post}(\nu_1,\nu_2 | \{q^i\}_{i=1}^{N_s})$ center more at this point.

The above approach can be extended if practitioner believes that the true pdf of $Z$ belongs to multiple families of distributions, each with its own set of parameters, i.e. $f^1_Z(z|\theta^1),\dots,f^{N_m}_Z(z| \theta^{N_m})$.  Bayesian model selection \cite{paper:GrigoriuF2008} offers a strategy to solve the inverse problem.

Since the practitioner possesses information about the family of distributions in which $Z$ resides, this approach guarantees that the true pdf of $Z$ can be recovered. Does this approach ensure that the optimization problem or the posterior density have a unique global minimum? As before, the answer is case dependent. For example, \cite{paper:LegollMOS2015} considers a specific model in stochastic homogenization where it was analytically shown that their specified distribution on $Z$ and the optimization problem they formulated accommodate a unique solution on the parameters of the pdf of $Z$.

Otherwise, additional information such as moments of $Z$ need to be supplied to regularize against other plausible parameter values in the pdf of $Z$. To clarify this, assume instead that the ranges of $Z_1,Z_2$ are unknown yet $f_Q(q) = -\log(q)$, $Z_1 = \lambda Z_1'$, $Z_2 = \frac{1}{\lambda} Z_2'$ where $\lambda > 0, Z_1' \sim Beta(\nu_1,\nu_1), Z_2' \sim Beta(\nu_2,\nu_2)$ are specified information. In this new model, $\lambda,\nu_1,\nu_2$ are parameters to be approximated.  Without additional information, the inverse problem possesses infinitely many solutions of the form $(\nu_1,\nu_2,\lambda) = (1,1,\lambda)$ for any $\lambda > 0$. Specifying additional moment information on $Z$ such as $E[Z_1]$ and $E[Z_2]$ would resolve such ill-posedness.

\section{Remarks}

Section~\ref{SISC:sec:AddlInfo} dealt with the objective of this work motivated by the discussion in Section~\ref{SISC:sec:AbsentInfo}. As we have seen, existing methods may not succeed in recovering the true pdf on the random quantity in the absence of further information. As such, we argued on how the inverse problem should be formulated and suggested appropriate solution methods. Here, we support and clarify aspects on the model construction (Section~\ref{SISC:subsec:ModelConstr}) and on our objective (Section~\ref{SISC:subsec:validation}). In particular, the latter section highlights that in the absence of information on the unknown random quantity, the solution to the inverse problem may not be suitable to predict the law of quantities of interest other than the one it was calibrated to.

\subsection{Posing the stochastic inverse problem} \label{SISC:subsec:ModelConstr}

In the following, three remarks are made about the various formulations of the inverse problem considered in Section~\ref{SISC:sec:AddlInfo}. The same notation above is used.


\begin{itemize}
\item {\it Different models on the pdf of $Z$ result in different solutions to the inverse problem}. To clarify this point, we revisit the examples presented in Sections~\ref{SISC:subsubsec:entropyBayes} and~\ref{SISC:subsubsec:KnownPdfBayes}. The parameters between both models are related through $\mu_i = \frac{1}{1 + \nu_i}$ and the same prior pdf characterized both sets of parameters. With the same set of $N_s = 100$ samples, we plot the posterior distributions under each model in the $(\mu_1,\mu_2)$ space (Figures~\ref{SISC:fig:EntropyPosterior} and~\ref{SISC:fig:KnownPdfPosteriorChangeVar}) and in the $(\nu_1,\nu_2)$ space  (Figures~\ref{SISC:fig:KnownPdfPosterior} and~\ref{SISC:fig:EntropyPosteriorChangeVar}). It is clear that the contours manifest distinct behavior.


\item {\it Difference in information content between knowing the pdf of $Q(Z)$ vs only having samples of $Q(Z)$.} Compared to having the pdf of $Q$ as given information, only possessing samples of $Q$ provides less information about the quantity of interest. A possible consequence of this lower information content includes posterior densities \eqref{SISC:eq:EntropyPosterior} that are not sharp about the MAP estimate. This is remedied by requiring a large number of samples of $Q(Z)$. 

\item {\it Distinction in applying the principle of maximum entropy when moment information on $Z$ or $Q(Z)$ is supplied.} Finally, we remark that the principle of maximum entropy has also been invoked to solve a stochastic inverse problem formulated differently from that of Section~\ref{SISC:sec:AddlInfo}. Consider the mapping $Q:D \times \Gamma \rightarrow \mathcal{D}$ with $\Gamma \coloneqq Z(\Omega)$, $D$ being the physical domain, and $Q(x,z) = q(U(x,z))$ for $x \in D, z \in \Gamma$ and some function $q$. The methodology developed in \cite{paper:GerbeauLT2018} seeks to address the inverse problem described as follows:

\begin{quote}
	Determine the pdf $f_Z$ of $Z$ given the bounded range $\Gamma$ of $Z$ and the observed moments of $Q$ up to order $N$ for $x\in D$ denoted by $\hat{\mu}^p(x)$, $p = 1,\dots,N$. 
\end{quote}

Since the $p$-th order moments of $Q$ can be expressed as integrals on $Z$, i.e. $ \displaystyle \int_{\Gamma} Q(x,z)^p f_Z(z) \,dz$ for $x \in D$, this naturally leads to a solution based on the principle of maximum entropy in which $f_Z$ is estimated via the optimization problem \eqref{SISC:eq:entropyProblem} subject to the constraint \\$\displaystyle \int_{\Gamma} Q(x,z)^p f_Z(z) \,dz = \hat{\mu}^p(x)   \,\,\,\forall x \in D, \,\,1 \le p \le N$.

Although the principle of maximum entropy has been used as a regularizer to infer the pdf of a random vector provided information about its moments, for problems involving forward models, the absence of information on $Z$ raises issues mentioned earlier. There is no guarantee that the proposed method is able to recover the true pdf of $Z$, in contrast to the authors' comment on p. B761. To see this, let $f^1_Z, f^2_Z$ be two pdfs of $Z$ such that when propagated through the model, the pdf of $Q$ is $f_Q$. The $p$-th order moments of $Q$ under both pdfs of $Z$ are identical even though $f^1_Z$ has larger entropy while $f^2_Z$ is the true pdf of $Z$ or vice versa. Specifying all moments of $Q$ neither resolves the issue.


\end{itemize}

\subsection{Validation} \label{SISC:subsec:validation}

In Section~\ref{SISC:sec:AddlInfo}, we considered types of information required on the unknown random quantity in order to solve the stochastic inverse problem. While some of this required information may be exigent, we argue that they are necessary to obtain solutions such that the resulting law of $Z$ can be used to characterize other quantities of interest $\widetilde{Q}$. We remark that in relation to the methods tackled in Section~\ref{SISC:sec:AbsentInfo}, such additional information may not be necessary if the structure of the contours of $Q$ and $\widetilde{Q}$ are similar. If this is not the case, without additional information, methods such as in Section~\ref{SISC:sec:ConsistentBayesian} may result in a posterior pdf for $Z$ whose predicted probability measure on the new quantity of interest $\widetilde{Q}$ is similar to the predicted measure on $\widetilde{Q}$ produced by the prior. The field of optimal experimental design for prediction addresses these concerns.

Here, we revisit the solutions obtained from the methods described in Sections~\ref{SISC:subsubsec:DesignLebesgueMeas}, \ref{SISC:subsubsec:entropyBayes}, and~\ref{SISC:subsubsec:KnownPdfBayes}. It is demonstrated that in the absence of information on $Z$, the resulting solution may be inadequate to characterize the law of quantities of interest to which it was not calibrated, thereby limiting its use in practical applications.

\begin{example} 
We revisit the forward mapping $Q(Z_1,Z_2) = Z_1 \cdot Z_2$ where $Z_1,Z_2 \sim U(0,1)$, independent, characterizes the true law on $Z$. The following methodologies are employed to solve the inverse problem on approximating the pdf of $Z$ depending on available information on $Z$ and $Q$. The resulting pdf on $Z$ is then used to predict the pdf on an unobserved quantity of interest $\widetilde{Q}(Z_1,Z_2) = Z_1 + Z_2$. The domain $Z \in [0,1]^2$ is assumed for all methods.

\begin{itemize}\label{SISC:ex:CompareSumQOI}
	\item Method based on the disintegration theorem using an ansatz as in Section~\ref{SISC:subsubsec:DesignLebesgueMeas} in which $f_Q(q) = -\log(q)$ is given and no other information on $Z$ is required.
	\item Bayes' theorem with entropy-based pdf for the likelihood as in Section~\ref{SISC:subsubsec:entropyBayes} in which $N_s = 100$ samples $\{q^i\}_{i=1}^{N_s}$ of $Q$ are available and the following is known about $Z$: $Z_1,Z_2$ are independent and $f_{\mu}^{prior}(\mu_1,\mu_2) = \mathbbm{1}_{(\mu_1,\mu_2) \in [0.25,0.75]^2} \frac{1}{0.5^2}$ is the prior pdf on $(\mu_1,\mu_2) = (E[Z_1],E[Z_2])$.
	\item Bayes' theorem with known family of distributions for the likelihood as in Section~\ref{SISC:subsubsec:KnownPdfBayes} in which the same $N_s = 100$ samples $\{q^i\}_{i=1}^{N_s}$ of $Q$ are available as above and the following is known about $Z$: $Z_1,Z_2$ are independent, $Z_1 \sim Beta(1,\nu_1), Z_2 \sim Beta(1,\nu_2)$, and $f_{\nu}^{prior}(\nu_1,\nu_2) = \mathbbm{1}_{(\nu_1,\nu_2) \in [0.75,1.25]^2}\frac{1}{0.5^2}$ is the prior pdf on $(\nu_1,\nu_2)$.
\end{itemize}
\end{example}

The solution approach for the latter 2 methods has already been discussed. The pdf on the unobserved quantity of interest $\widetilde{Q}$ then results by computing the posterior predictive distribution. Denote by $f_{\Theta}^{post}(\theta | \{q^i\}_{i=1}^{N_s})$ the obtained posterior pdf on the corresponding parameter space which qualifies as the solution to the inverse problem upon application of either of the latter 2 methods. The pdf on $\widetilde{Q}(Z_1,Z_2) = Z_1+Z_2$ is obtained through
\begin{align}\label{SISC:eq:posteriorPredictiveDist}
f_{\widetilde{Q}}(\widetilde{q}| \{q^i\}_{i=1}^{N_s}) = \int_{\Theta} f_{\widetilde{Q}}(\widetilde{q}| \theta) \cdot f_{\Theta}^{post}(\theta | \{q^i\}_{i=1}^{N_s}) \,d\theta 
\end{align}
where $f_{\widetilde{Q}}(\widetilde{q}|\theta)$ is the pdf on $\widetilde{Q}$ obtained by propagating the conditional pdf $f_Z(\cdot| \theta)$ on $Z$ through $\widetilde{Q}$. Elementary calculations show that
\begin{equation} \label{SISC:eq:SumQOIPdf}
f_{\widetilde{Q}}(\widetilde{q}|\theta) = \left\{
\begin{array}{ll} 
      \displaystyle \int_0^{\widetilde{q}} f_{Z_2}(\widetilde{q} - z_1 | \theta) \, f_{Z_1}(z_1| \theta) \,dz_1 & 0 \le \widetilde{q} \le 1 \\
      \displaystyle \int_{\widetilde{q}-1}^1 f_{Z_2}(\widetilde{q} - z_1 | \theta) \, f_{Z_1}(z_1| \theta) \,dz_1 & 1 \le \widetilde{q} \le 2 \\
\end{array} 
\right.
\end{equation}
with $f_{Z_i}(\cdot|\theta)$ being the marginal pdf of $Z_i$, $i=1,2$. In particular, under the true law on $Z$, $\widetilde{Q}$ has a triangular distribution whose pdf is  $f_{\widetilde{Q}}(\widetilde{q}) = \widetilde{q}$ for $\widetilde{q} \in [0,1]$ and $f_{\widetilde{Q}}(\widetilde{q}) = 2-\widetilde{q}$ for $\widetilde{q} \in [1,2]$.

On the other hand, the pdf $f_Z^{ansatz}$ arising from the method based on the disintegration theorem results through this procedure:
\begin{enumerate}
	\item Partition $\Gamma = [0,1]^2$ into $N_{sq}=10000$ squares $\{A_i\}_{i=1}^{N_{sq}}$ with area $(0.01)^2$.
	\item Compute $P_Z(A_i)$ using \eqref{SISC:eq:disIntForProbLebesgue} and \eqref{SISC:eq:ansatzLebesgue} for $i=1,\dots,N_{sq}$.
	\item For each $A_i$, let $(z_{1,i}^*,z_{2,i}^*)$ be its center. The approximate pdf of $Z$ is then calculated as $f_Z^{ansatz}(z_{1,i}^*,z_{2,i}^*) \simeq \frac{P_Z(A_i)}{(0.01)^2}$.
\end{enumerate}
	
If $(z_{1,i}^{SW},z_{2,i}^{SW})$ and $(z_{1,i}^{NE},z_{2,i}^{NE})$ represent the lower left and the upper right vertices, respectively, of $A_i$, the parameter along the transverse curve $x_{\mathcal{L}}$ is bounded by $ \sqrt{2z_{1,i}^{SW}z_{2,i}^{SW}}\le x_{\mathcal{L}} \le \sqrt{2z_{1,i}^{NE}z_{2,i}^{NE}}$.

\begin{figure}[h!]
		\centering
		\includegraphics[width = 0.9\textwidth]		
						{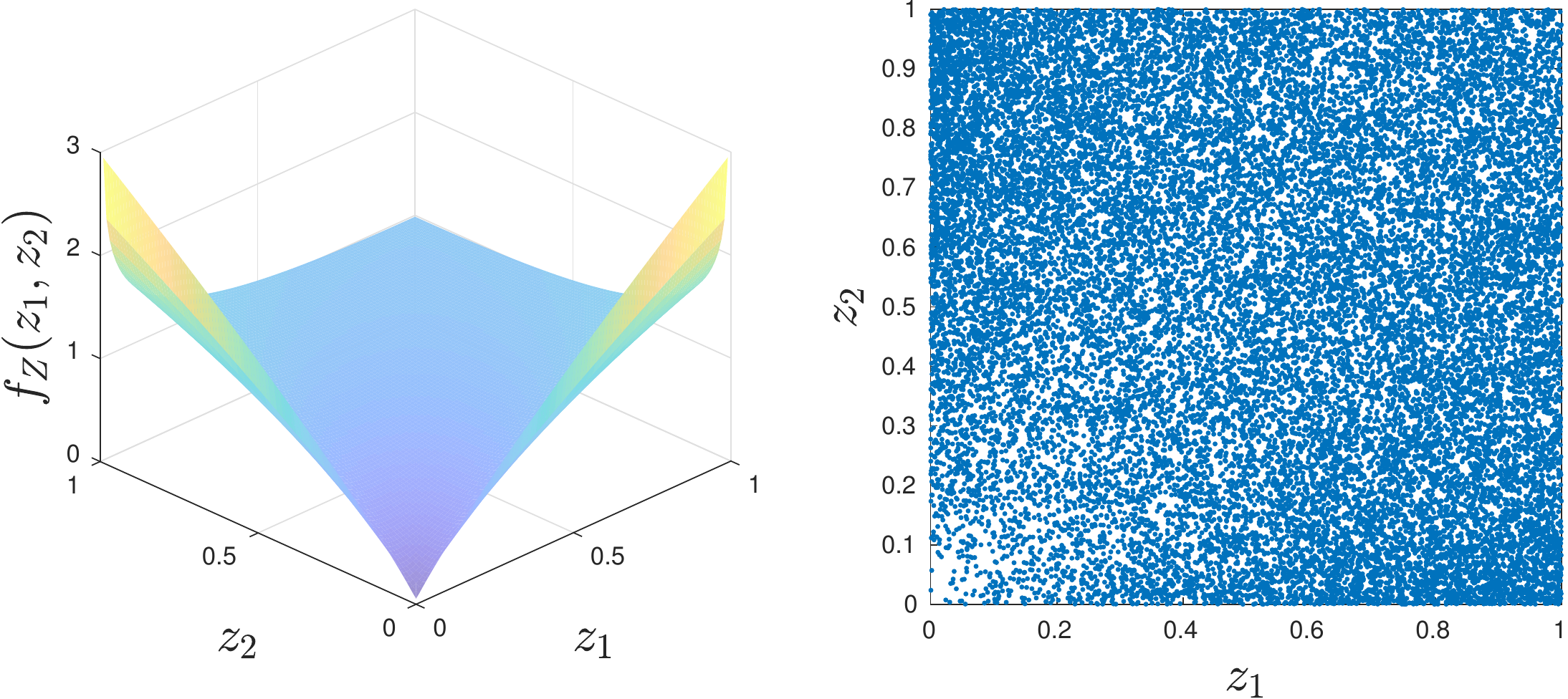}
		\caption{Left panel: approximate pdf of $Z$ produced by the method in Section~\ref{SISC:subsubsec:DesignLebesgueMeas}. Right panel: 25000 samples of $(Z_1,Z_2)$ simulated from the pdf on the left panel.}
		\label{SISC:fig:BETMethodAndSamples}
\end{figure}

\begin{figure}[h!]
		\centering
		\includegraphics[width = 0.4\textwidth]		
						{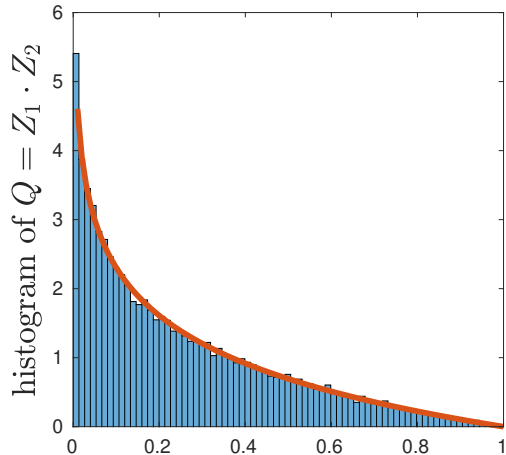}
		\caption{Histogram of 52154 samples of $Q = Z_1\cdot Z_2$ where the samples of $Z$ are drawn from $f_Z^{ansatz}$ together with $f_Q(q) = -\log(q)$.}
		\label{SISC:fig:ObsQOIRejSampling}
\end{figure}

The left panel of Figure~\ref{SISC:fig:BETMethodAndSamples} displays the approximate pdf of $Z$ produced by the method in Section~\ref{SISC:subsubsec:DesignLebesgueMeas} whereas the right panel shows 25000 samples of $(Z_1,Z_2)$ drawn from this pdf through rejection sampling and interpolation. We confirm that when $f_Z^{ansatz}$ is propagated through the forward model $Q$, we recover the specified pdf $f_Q(q) = -\log(q)$ as guaranteed by Theorem~\ref{SISC:thm:DisintThm}. Figure~\ref{SISC:fig:ObsQOIRejSampling} exhibits the histogram of $Z_1\cdot Z_2$ using 52154 samples of $Z$ obtained from $f_Z^{ansatz}$ together with $f_Q$.


With the pdf $f_Z^{ansatz}$ on $Z$ at hand, one approach to propagate this pdf through the unobserved quantity of interest $\widetilde{Q}$ is to construct a discrete pdf approximation to $\widetilde{Q}$ using the centers of each square $(z_{1,i}^*,z_{2,i}^*)$ and their corresponding probabilities $P_Z(A_i)$. This yields distinct outcomes $\widetilde{q}_j$ of $\widetilde{Q}$ with weights $P(\widetilde{Q} = \widetilde{q}_j) = \sum_{i: z_{1,i}^*+z_{2,i}^* = \widetilde{q}_j} P_Z(A_i)$ that need to normalized to obtain the approximate pdf $f_{\widetilde{Q}}^{ansatz}$ of $\widetilde{Q}$.

\begin{figure}[h!]
		\centering
		\includegraphics[width = 1.05\textwidth]		
						{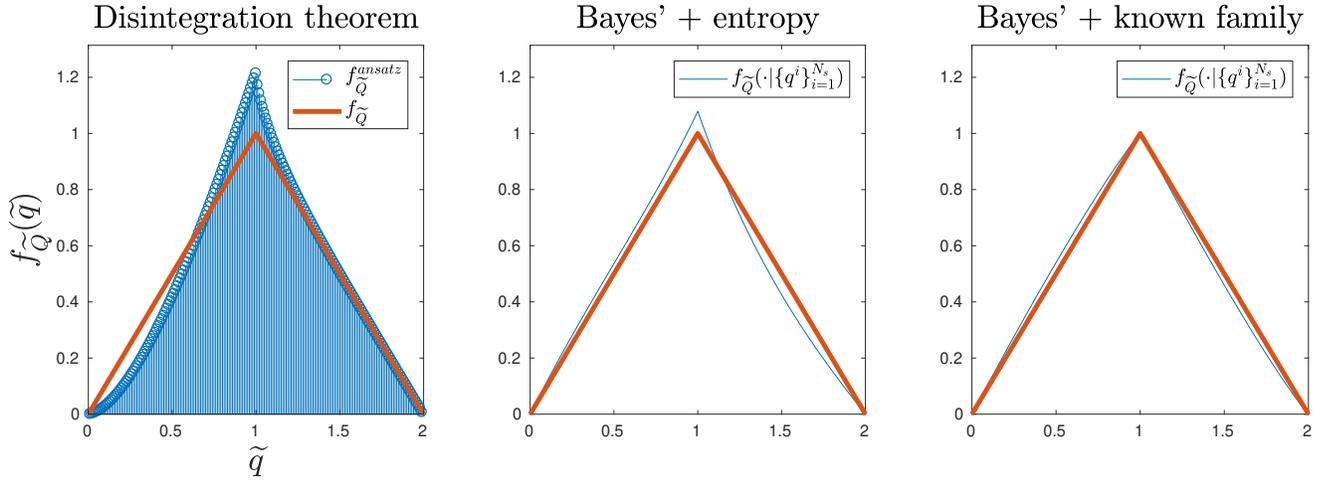}
		\caption{Comparison between the true pdf $f_{\widetilde{Q}}$ of $\widetilde{Q}$ and the pdf obtained by propagating through $\widetilde{Q}$ the pdf of $Z$ stemming from methods for the inverse problem described in Example~\ref{SISC:ex:CompareSumQOI}. A stem plot is used for $f_{\widetilde{Q}}^{ansatz}$ in the left-most plot to emphasize that an accurate discrete random variable approximation was used.}
		\label{SISC:fig:CompareMethodsQOISum}
\end{figure}

We now evaluate the performance of the solution from each method to predict the probability law of the unobserved $\widetilde{Q}$. Figure~\ref{SISC:fig:CompareMethodsQOISum} contains 3 subplots, one for each method, plotting $f_{\widetilde{Q}}^{ansatz}$ or $f_{\widetilde{Q}}(\cdot|\{q^i\}_{i=1}^{N_s})$ together with the true pdf of $\widetilde{Q}$. A stem plot was used for the leftmost subplot to emphasize that the pdf of a discrete random variable was used to accurately approximate $f_{\widetilde{Q}}^{ansatz}$. Quantitatively, the discrepancy between the true pdf and the simulated pdf's is: $\text{max} |f_{\widetilde{Q}}^{ansatz} - f_{\widetilde{Q}}| \approx 0.2141$ whereas $\text{max} |f_{\widetilde{Q}}(\cdot|\{q^i\}_{i=1}^{N_s}) - f_{\widetilde{Q}}| \approx 0.0786$ for Method 2 (Bayes'+entropy) while $\text{max} |f_{\widetilde{Q}}(\cdot|\{q^i\}_{i=1}^{N_s}) - f_{\widetilde{Q}}| \approx 0.0415$ for Method 3 (Bayes'+known family). It was also observed that increasing the number of available samples $\{q^i\}_{i=1}^{N_s}$ of the observed quantity of interest $Q$ decreased the discrepancy for the latter 2 methods. The objective of this example was not to conclude which method is better since the available information for each was not identical. Rather, this example justifies the need to specify additional information on $Z$ as was carried out in Sections~\ref{SISC:sec:AddlInfo}.

We conclude this section by showing the repercussions that may arise if the pdf $f_Z^{ansatz}$ on $Z$ is used to predict the pdf of more complicated quantities of interest $\widetilde{Q}$. If we consider $\widetilde{Q}(Z_1,Z_2) = Z_1^2 + Z_2^2$ where $Z_1,Z_2 \sim U(0,1)$ and independent, it was proven in \cite{paper:Weissman2017} that 

\begin{equation*}
f_{\widetilde{Q}}(\widetilde{q}) = \left\{
\begin{array}{ll} 
      \frac{\pi}{4} & 0 \le \widetilde{q} \le 1 \\
      \arcsin \frac{1}{\sqrt{\widetilde{q}}} - \frac{\pi}{4} & 1 \le \widetilde{q} \le 2 \\
\end{array} 
\right..
\end{equation*}

\begin{figure}[h!]
		\centering
		\includegraphics[width = 1\textwidth]		
						{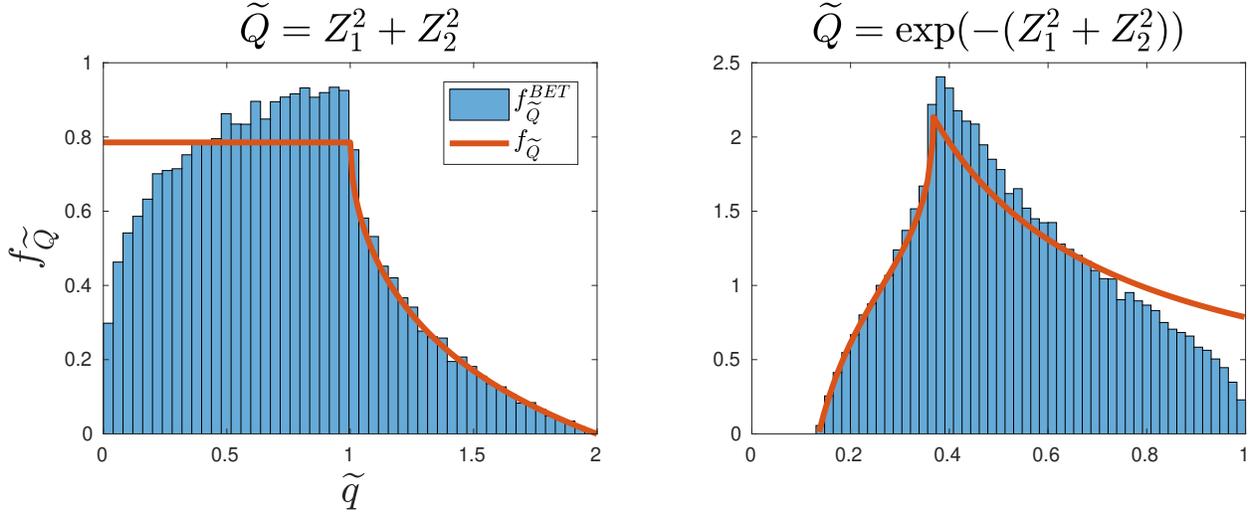}
		\caption{Comparison between the true pdf $f_{\widetilde{Q}}$ of $\widetilde{Q}$ and $f_{\widetilde{Q}}^{ansatz}$ for more complicated quantities of interest.}
		\label{SISC:fig:MoreComplexQOI}
\end{figure}

Additionally, we also consider $\widetilde{Q}(Z_1,Z_2) = \exp(-(Z_1^2 + Z_2^2))$. Figure~\ref{SISC:fig:MoreComplexQOI} contains subplots comparing the true pdf $f_{\widetilde{Q}}$ and histogram of $f_{\widetilde{Q}}^{ansatz}$ computed via samples of $Z$ drawn from $f_{Z}^{ansatz}$ for the two aforementioned unobserved $\widetilde{Q}$'s. We see that $f_{\widetilde{Q}}^{ansatz}$ under/overestimates probabilistic properties of $\widetilde{Q}$ such as $P(\widetilde{Q} \le 0.5)$ for the left subplot and $P(\widetilde{Q} \ge 0.8)$ for the right subplot, among other properties such as moments of $\widetilde{Q}$. 

The above examples are not intended to discredit existing methods for solving stochastic inverse problems. Instead, they underscore why additional information should be specified on $Z$ to yield solutions that can be useful in applications. In fact, existing works such as \cite{paper:NolenP2009,paper:LegollMOS2015,paper:EmeryGF2016} impose probability distributions on the unknown random quantities after which standard mathematical tools were applied to infer the  parameters that characterize these distributions.

\section{Conclusion}

This work dealt with the stochastic inverse problem of identifying the distribution of $Z$ given probabilistic information of the quantity of interest $Q(Z)$. We surveyed general methods that have been developed to tackle the problem in which no information other than the bounded range of $Z$ is assumed. These methods coped with the ill-posedness of this inverse problem in the following ways: \cite{paper:BreidtBE2011, paper:ButlerETDW2014, paper:ButlerGEDW2015} suggested an ansatz that the pdf on the generalized contours of $Q$ is uniform while \cite{paper:BorggaardV2015} obtained the probability law of $Z$ from the Karhunen-Lo\`eve expansion of the solution to the stochastic PDE. We motivated this work by showing that this lack of additional information entails that the true pdf of $Z$ may not be recovered. 

Consequently, we argued that it is necessary for this inverse problem to be posed such that further information on $Z$ is specified to attain solutions that are of practical use. We demonstrated that this specified information can take the form of moments of $Z$ or the family of distributions in which $Z$ resides subject to unknown parameters, among others. Using these information, a conjunction of tools such as Bayes' theorem, the principle of maximum entropy, and forward uncertainty propagation were utilized to solve the inverse problem in a manner that is consistent with the present information on the quantity of interest and on $Z$. Issues arising from this framework were also highlighted. Finally, we emphasized the need for this specified information by assessing how well the resulting solutions from the discussed methods can predict the probability law of unobserved quantities of interest. It is observed that solving the inverse problem without additional information on $Z$ can lead to solutions that may be unreliable for prediction.

The intention of this work was not to discredit existing contributions in this area but to stress on how we believe the inverse problem must be posed for use in practical applications.

\appendix 
\numberwithin{equation}{section}

\section{Computing the pdf on 1-dimensional contours for Example 1} \label{SISC:appendix:pdfContour}

Fix $x_{\mathcal{L}}$ and consider the contour $\pi^{-1}(x_{\mathcal{L}})$. Let $x_\mathcal{C}$ parameterize the arc length of  $\pi^{-1}(x_{\mathcal{L}})$ and denote by $\mu(\pi^{-1}(x_{\mathcal{L}}))$ the arc length of $\pi^{-1}(x_{\mathcal{L}})$ so that $0 \le x_\mathcal{C} \le \mu(\pi^{-1}(x_{\mathcal{L}}))$. The conditional pdf
$f_{X_{\mathcal{C}}|X_{\mathcal{L}}} (x_{\mathcal{C}}|x_{\mathcal{L}})$ along the contour can be approximated as follows:

\begin{itemize}
	\item Select equidistant points $\{x_{\mathcal{C}}^{(i)}\}_{i=0}^N \subset [0,\mu(\pi^{-1}(x_{\mathcal{L}}))] $ such that $x_{\mathcal{C}}^{(0)} = 0, x_{\mathcal{C}}^{(N)} = \mu(\pi^{-1}(x_{\mathcal{L}}))$ and  $x_{\mathcal{C}}^{(i+1)} - x_{\mathcal{C}}^{(i)} = \frac{\mu(\pi^{-1}(x_{\mathcal{L}}))}{N}$ for $i = 0,\dots,N-1$.
	\item From $P(X_{\mathcal{C}} \in (x_{\mathcal{C}}^{(i)},x_{\mathcal{C}}^{(i+1)}) \, | \, X_{\mathcal{L}} = x_{\mathcal{L}}) = \displaystyle \int_{x_{\mathcal{C}}^{(i)}}^{x_{\mathcal{C}}^{(i+1)}} f_{X_{\mathcal{C}} \,| \,X_{\mathcal{L}}} (x_{\mathcal{C}}|x_{\mathcal{L}})\,  dx_{\mathcal{C}}$, we deduce the approximation
\begin{align} \label{SISC:eq:contourPDFApprox}
f_{X_{\mathcal{C}}|X_{\mathcal{L}}} (x_{\mathcal{C}}|x_{\mathcal{L}}) \simeq \frac{P(X_{\mathcal{C}} \in (x_{\mathcal{C}}^{(i)},x_{\mathcal{C}}^{(i+1)}) \, | \, X_{\mathcal{L}} = x_{\mathcal{L}})}{x_{\mathcal{C}}^{(i+1)}-x_{\mathcal{C}}^{(i)}}
\end{align}
for $x_{\mathcal{C}} \in (x_{\mathcal{C}}^{(i)},x_{\mathcal{C}}^{(i+1)})$ assuming that $N$ is sufficiently large.
\item To approximate the numerator in \eqref{SISC:eq:contourPDFApprox}, we construct an infinitesimal region $R$ in $\Gamma = [0,1]^2$ bounded by the contours $\pi^{-1}(x_{\mathcal{L}})$ and $\pi^{-1}(x_{\mathcal{L}}+\epsilon)$ for $\epsilon$ sufficiently small and partition $R$ into regions $\{R_i\}_{i=1}^N \subset \Gamma$ with $R_i$ corresponding to the arc lengths between $x_{\mathcal{C}}^{(i-1)}$ and $x_{\mathcal{C}}^{(i)}$ for $i=1,\dots,N$. It then follows that
\begin{align} \label{SISC:eq:integralExpression}
P(X_{\mathcal{C}} \in (x_{\mathcal{C}}^{(i)},x_{\mathcal{C}}^{(i+1)}) \, | \, X_{\mathcal{L}} = x_{\mathcal{L}}) \simeq \displaystyle \iint_{R_{i+1}} f_Z(z_1,z_2) \,\, dz_1 dz_2
\end{align}
where $f_Z$ is the joint pdf of $Z_1,Z_2$.
\end{itemize}


\begin{thebibliography}{9}

\bibitem{book:Grigoriu2012}
M. Grigoriu, \emph{Stochastic systems. Uncertainty quantification and propagation}, Springer Ser. Reliab. Eng., Springer, London, 2012.

\bibitem{book:SoongG1993}
T.T. Soong and M. Grigoriu, \emph{Random vibrations of structural and mechanical systems}, Prentice Hall, New Jersey, 1993.

\bibitem{book:Grigoriu2002}
M. Grigoriu, \emph{Stochastic calculus. Applications in science and engineering.}, Birkh\"auser, Boston, 2002.

\bibitem{book:Hansen2010}
P.C. Hansen, \emph{Discrete Inverse Problems: Insight and Algorithms}, SIAM, Philadelphia, 2010.

\bibitem{book:KaipioS2005}
J. Kaipio and E. Somersalo, \emph{Statistical and Computational Inverse Problems}, Springer, New York, 2005.

\bibitem{book:WongH1985}
E. Wong and B. Hajek, \emph{Stochastic processes in engineering systems}, Springer, New York, 1985.

\bibitem{book:CoverT2006}
T.M. Cover and J.A. Thomas, \emph{Elements of information theory}, 2nd ed., Wiley, New York, 2006.

\bibitem{book:Durrett2010}
R. Durrett, \emph{Probability: Theory and examples}, 4th ed., Cambridge University Press, New York, 2010.

\bibitem{paper:BreidtBE2011}
J. Breidt, T. Butler, and D. Estep, \emph{A measure-theoretic computational method for inverse sensitivity problems I: method and analysis}, SIAM J. Numerical Analysis 49 (2011) 1836-1859.

\bibitem{paper:ButlerETDW2014}
T. Butler, D. Estep, S. Tavener, C. Dawson, and J. J. Westerink, \emph{A measure-theoretic computational method for inverse sensitivity problems III: multiple quantities of interest}, SIAM/ASA J. Uncertainty Quantification 2 (2014) 174-202.

\bibitem{paper:ButlerGEDW2015}
T. Butler, L. Graham, D. Estep, C. Dawson, and J. J. Westerink, \emph{Definition and solution of a stochastic inverse problem for the Manning's \it{n} parameter field in hydrodynamic models}, Advances in Water Resources 78 (2015) 60-79.

\bibitem{paper:MattisBDEV2015}
S.A. Mattis, T.D. Butler, C.N. Dawson, D. Estep, and V.V. Vesselinov, \emph{Parameter estimation and prediction for groundwater contamination based on measure theory}, Water Resour Res 51(9) (2015) 7608-7628.

\bibitem{paper:FieldG2007}
R. V. Field and M. Grigoriu, \emph{Convergence Properties of Polynomial Chaos Approximations for $L^2$-Rndom Variables}. Sandia Report SAND2007-1262, 2007.

\bibitem{paper:FieldG2005}
R.V. Field Jr. and M. Grigoriu, \emph{On the Accuracy of the Polynomial Chaos Approximation}, J. Comput. Phys. 209 (2005) 617-642.

\bibitem{paper:BorggaardV2015}
J. Borggaard and H. van Wyk, \emph{Gradient-based estimation of uncertain parameters for elliptic partial differential equations}, Inverse Problems 31 (2015) 065008.

\bibitem{paper:BanksB2001}
H.T. Banks and K.L. Bihari, \emph{Modelling and estimating uncertainty in parameter estimation}, Inverse Problems 17 (2001) 95-111.

\bibitem{paper:DesceliersGS2006}
C. Desceliers, R. Ghanem, C. Soize, \emph{Maximum likelihood estimation of stochastic chaos representations from experimental data}, Int. J. Numer. Meth. Engng 66 (2006) 978-1001.

\bibitem{paper:ZabarasG2008}
N. Zabaras and B. Ganapathysubramanian, \emph{A scalable framework for the solution of stochastic inverse problems using a sparse grid collocation approach}, Journal of Computational Physics 227 (2008) 4697-4735.

\bibitem{paper:NarayananZ2004}
V.A.B. Narayanan and N. Zabaras, \emph{Stochastic inverse heat conduction using a spectral approach}, Int. J. Numer. Meth. Engng 60 (2004) 1-24.

\bibitem{paper:WarnerAG2014}
J.E. Warner, W. Aquino, and M.D. Grigoriu, \emph{Stochastic reduced order models for inverse problems under uncertainty}, Comput. Methods Appl. Mech. Engrg. 285 (2015) 488-514.

\bibitem{paper:BabuskaTZ2005}
I. Babu\v{s}ka, R. Tempone, and G.E. Zouraris, \emph{Solving elliptic boundary value problems with uncertain coefficients by the finite element method: the stochastic formulation}, Comput. Methods in Appl. Mech. Eng. 194 (2005) 1251-1294. 

\bibitem{paper:GerbeauLT2018}
J. Gerbeau, D. Lombardi, and E. Tixier, \emph{A moment-matching method to study the variability of phenomena described by partial differential equations}, SIAM J. Sci. Comput. 40(3) (2018) B743-B765.

\bibitem{paper:ButlerJW2018}
T. Butler, J. Jakeman, and T. Wildey, \emph{Combining push-forward measures and Bayes' rule to construct consistent solutions to stochastic inverse problems}, SIAM J. Sci. Comput. 40(2) (2018) A984-A1011.

\bibitem{paper:NolenP2009}
J. Nolen and G. Papanicolaou, \emph{Fine scale uncertainty in parameter estimation for elliptic equations}, Inverse Problems 25 (2009) 115021.

\bibitem{paper:LegollMOS2015}
F. Legoll, W. Minvielle, A. Obliger, and M. Simon, \emph{A parameter identification problem in stochastic homogenization}, ESAIM Proc. 48 (2015) 190-214.

\bibitem{paper:EmeryGF2016}
J.M. Emery, M.D. Grigoriu, R.V. Field Jr, \emph{Bayesian methods for characterizing unknown parameters of material models}, Applied Mathematical Modelling 40(13) (2016) 6395-6411.

\bibitem{paper:Weissman2017}
I. Weissman, \emph{Sum of squares of uniform random variables}, Statistics and Probability Letters 129 (2017) 147-154.

\bibitem{paper:GrigoriuF2008}
M.D. Grigoriu and R.V. Field Jr, \emph{A solution to the static frame validation challenge problem using Bayesian model selection}, Comput. Methods. Appl. Mech. Engrg. 197 (2008) 2540-2549.

\end{thebibliography}
\end{document}